\documentclass[12pt, reqno, microtype]{amsart}

\usepackage{amsmath, amssymb, amsthm}
\usepackage[hidelinks, unicode]{hyperref}
\usepackage[parfill]{parskip}
\usepackage{pinlabel}
\usepackage[all]{xy}
\usepackage[margin=1.05in]{geometry}

\usepackage{tikz}
\usetikzlibrary{matrix,arrows,positioning}

\usepackage{verbatim} 

\usepackage{mathrsfs}
\usepackage{mathtools}
\usepackage{mathabx}

\title{Pattern-equivariant homology}

\author{James J.\ Walton} 
\address{Department of Mathematics\\University of York\\\newline
         Heslington, York, YO10 5DD\\UK}
\email{jamie.walton@york.ac.uk}
\urladdr{http://maths.york.ac.uk/www/jjw548}

\thanks{Research supported by EPSRC}

\theoremstyle{definition}

\numberwithin{equation}{section}

\newcommand{\mf}{\mathfrak}
\newcommand{\mc}{\mathcal}

\newcommand{\mathbbm}[1]{\text{\usefont{U}{bbm}{m}{n}#1}}

\newcommand{\R}{\mathbb{R}}
\newcommand{\Q}{\mathbb{Q}}
\newcommand{\Z}{\mathbb{Z}}
\newcommand{\N}{\mathbb{N}}
\newcommand{\cf}{\mathbbm{1}}

\def\co{\colon\thinspace}

\def\co{\colon\thinspace}

\newtheorem{theorem}{Theorem}[section]
\newtheorem*{theorem*}{Theorem}
\newtheorem{proposition}[theorem]{Proposition}
\newtheorem{lemma}[theorem]{Lemma}
\newtheorem*{lemma*}{Lemma}

\newtheorem*{corollary*}{Corollary}

\theoremstyle{definition}

\newtheorem{definition}[theorem]{Definition}
\newtheorem*{definition*}{Definition}
\newtheorem{exmp}[theorem]{Example}

\numberwithin{equation}{section} 
\numberwithin{figure}{section}   

\graphicspath{ {images/} }

\begin{document}
\begin{abstract}
Pattern-equivariant (PE) cohomology is a well-established tool with which to interpret the \v{C}ech cohomology groups of a tiling space in a highly geometric way. In this paper we consider homology groups of PE infinite chains. We establish Poincar\'{e} duality between the PE cohomology and PE homology. The Penrose kite and dart tilings are taken as our central running example, we show how through this formalism one may give highly approachable geometric descriptions of the generators of the \v{C}ech cohomology of their tiling space. These invariants are also considered in the context of rotational symmetry. Poincar\'{e} duality fails over integer coefficients for the `ePE homology groups' based upon chains which are PE with respect to orientation-preserving Euclidean motions between patches. As a result we construct a new invariant, which is of relevance to the cohomology of rotational tiling spaces. We present an efficient method of computation of the PE and ePE (co)homology groups for hierarchical tilings.
\end{abstract}
\maketitle

\section*{Introduction}

In the past few decades a rich class of highly ordered patterns has emerged whose central examples, despite lacking global translational symmetries, exhibit intricate internal structure, imbuing these patterns with properties akin to those enjoyed by periodically repeating patterns. The field of \emph{aperiodic order} aims to study such patterns, and to establish connections between their properties, and their constructions, to other fields of mathematics and the natural sciences. To name a few, aperiodic order has interactions with areas of mathematics such as mathematical logic \cite{LafWei08}---as established by Berger's proof of the undecidability of the domino problem \cite{Berg66}, Diophantine approximation \cite{ABEI01, BerSie05, HKW14, HKW15}, the structure of attractors \cite{ClaHun12} and symbolic dynamics \cite{Sch01}, and notably is of relevance to solid state physics in the wake of the discovery of quasicrystals by Shechtman et al.\ \cite{SBGC84}.

A full understanding of a periodic tiling, modulo locally defined reversible redecorations, amounts to an understanding of its symmetry group. In the aperiodic setting, the complexity and incredible diversity of examples demands a multifaceted approach. Techniques from the theory of groupoids \cite{BJS10}, semigroups \cite{KelLaw00}, $C^*$-algebras \cite{AndPut98}, dynamical systems \cite{ClaSad06, Kel08}, ergodic theory \cite{Rad97} and shape theory \cite{ClaHun12} find natural r\^{o}les in the field, and of course these tools have tightly knit connections to each other \cite{KelPut00}. One approach to studying a given aperiodic tiling $\mf{T}$ is to associate to it a moduli space $\Omega$, sometimes called the \emph{tiling space}, of locally indistinguishable tilings imbued with a natural topology; see Sadun's book \cite{Sad08} for an accessible introduction to the theory. A central goal is then to formulate methods of computing topological invariants of $\Omega$, and to describe what these invariants actually tell us about the original tiling $\mf{T}$. An important perspective, particularly for the latter half of this objective, is provided by Kellendonk and Putnam's theory of pattern-equivariant (PE) cohomology \cite{Kel03,KelPut06}. PE cohomology allows for an intuitive description of the \v{C}ech cohomology $\check{H}^\bullet(\Omega)$ of tiling spaces. Over $\R$ coefficients the PE cochain groups may be defined using PE differential forms \cite{Kel03}, and over general Abelian coefficients, when the tiling has a cellular structure, with PE cellular cochains \cite{Sad07}. Rather than just providing a reflection of topological invariants of tiling spaces, on the contrary, these PE invariants are of principle relevance to aperiodic structures and their connections with other fields in their own right; see, for example, Kelly and Sadun's use of them in a topological proof of theorems of Kesten and Oren regarding the discrepancy of irrational rotations \cite{KelSad15}. It is perhaps more appropriate to view the isomorphism between $\check{H}^\bullet(\Omega)$ and the PE cohomology as an elegant interpretation of the \emph{PE cohomology}, rather than vice versa, of theoretical and computational importance.

\begin{figure}
\begin{center}
\includegraphics[width=\textwidth]{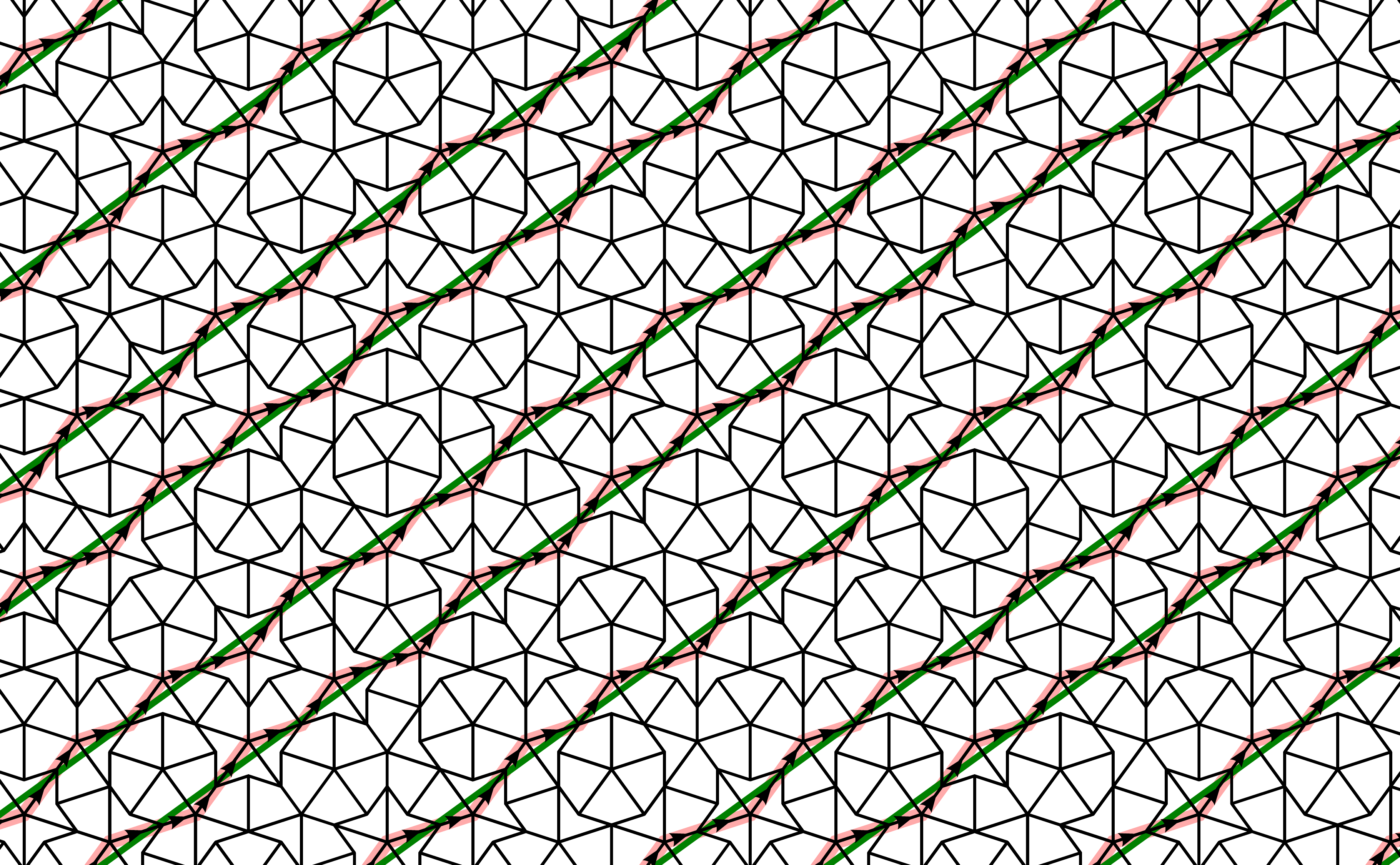}
\caption{The PE $1$-cycle $\nu_0$, with green Ammann bars of the supertiling.}
\label{fig: Ammann Chain}
\end{center}
\end{figure}

In this paper we introduce the \emph{pattern-equivariant homology groups} of a tiling. These homology groups are based on infinite, or non-compactly supported cellular chains, sometimes known as `Borel--Moore chains'. We say that such a chain is \emph{pattern-equivariant} if there exists some $r>0$ for which the coefficient of a cell only depends on the translation class of that cell and its surrounding patch of tiles to radius $r$. We show in Theorem \ref{thm: PE PD}, via a classical `cell, dual-cell' argument, that for a tiling of finite local complexity (see Subsection \ref{subsec: Cellular, Polytopal and Dual Tilings}) we have PE Poincar\'{e} duality:

{\bf Theorem \ref{thm: PE PD}}
{\sl For a polytopal tiling $\mf{T}$ of $\R^d$ of finite local complexity, we have PE Poincar\'{e} duality $H^\bullet(\mf{T}^1) \cong H_{d-\bullet}(\mf{T}^1)$ between the PE cohomology and PE homology of $\mf{T}$.}

The upshot of this is that one may give quite beautiful, and informative, geometric depictions of the \v{C}ech cohomology groups $\check{H}^\bullet(\Omega)$ of tiling spaces. For example, the cohomology of the translational hull $\Omega_\mf{T}^1$ of a Penrose kite and dart tiling in degree one is $\check{H}^1(\Omega_\mf{T}^1) \cong H^1(\mf{T}^1) \cong H_1(\mf{T}^1) \cong \Z^5$. The generators of this group have down-to-Earth interpretations in terms of important geometric features of the Penrose tilings. For example, one such generator, depicted in Figure \ref{fig: Ammann Chain}, is closely linked to Ammann bars of the Penrose tilings (of which, see the discussion in \cite[Chapter 10.6]{GruShe87}). Another simple geometric feature of the Penrose tilings is that the dart tiles arrange as loops, leading to the cycle depicted in Figure \ref{fig: Penrose}. As described in Example \ref{ex: Penrose PE} these two chains, and close analogues of them, give a near complete description of $H_1(\mf{T}^1)$.

In Section \ref{sec: PE Homology and Rotations} we consider these PE invariants in the context of rotational symmetry. Whilst for a tiling of finite local complexity the action of rotation on the PE homology and cohomology agree via Poincar\'{e} duality (Proposition \ref{prop: rotation actions}), the actions at the (co)chain level behave differently. We consider \emph{ePE} chains and cochains, which are required to have the same coefficients at any two cells whenever those cells agree on patches of sufficiently large radius up to orientation-preserving Euclidean motion (rather than just translations as in the case of the PE homology groups). We show in Theorem \ref{thm: ePE PD} that over divisible coefficients $G$ we still have Poincar\'{e} duality $H^\bullet(\mf{T}^0;G) \cong H_{d-\bullet}(\mf{T}^0;G)$ between the ePE cohomology and ePE homology, but over $\Z$ coefficients this typically fails. For example, for the Penrose kite and dart tilings the degree zero ePE homology group has a copy of an order five cyclic subgroup not present in the corresponding ePE cohomology group in degree two. A degree zero ePE torsion element is depicted in Figure \ref{fig: S0TorBdy}.

So whilst the PE homology gives a curious alternative way of visualising PE invariants, the ePE homology provides a new invariant to the ePE cohomology, or the \v{C}ech cohomology of the associated space $\Omega_\mf{T}^0$ (defined in \cite{BDHS10}, or see Subsection \ref{subsec: Rotational Tiling Spaces}). We shall show in forthcoming work \cite{Wal16} how this new invariant may naturally be incorporated into a spectral sequence converging to the \v{C}ech cohomology of the `Euclidean hull' $\Omega_\mf{T}^{\text{rot}}$ (see Subsection \ref{subsec: Rotational Tiling Spaces}) of a two-dimensional tiling. The only potentially non-trivial map of this spectral sequence has a very simple description in terms of the local combinatorics of the tiling. This procedure dovetails conveniently with the methods that we shall introduce in Section \ref{sec: PE Homology of Hierarchical Tilings} to efficiently compute the \v{C}ech cohomology of Euclidean hulls of hierarchical tilings, leading to some new computations on the cohomologies of these spaces.

We show how the ePE homology, ePE cohomology and rotationally invariant part of the PE cohomology are related for a two-dimensional tiling in Theorem \ref{thm: rotationally invariant triangle} and give the corresponding calculations for the Penrose kite and dart tilings. In general, over rational coefficients all three are canonically isomorphic, but over integral coefficients the canonical map from the ePE cohomology to the rotationally invariant part of the PE cohomology is rarely an isomorphism. It turns out that this map naturally factorises through the ePE homology. In some sense, the ePE homology adds extra cycles to the ePE cohomology.

The techniques that we present are not limited to tilings of Euclidean space. In Subsection \ref{subsec: Generalising} we introduce the notion of a \emph{system of internal symmetries}, which neatly encodes the necessary data required to define PE cohomology and various other related constructions. This allows us, for example, to apply the same techniques to non-Euclidean tilings, such as the combinatorial pentagonal tilings \cite{BowSte97} of Bowers and Stephenson.

In Section \ref{sec: PE Homology of Hierarchical Tilings} we change tack by considering the problem of how to actually compute the PE homology for certain examples. The PE homology formalism naturally leads to a simple, and efficient method of computation for invariants of a hierarchical tilings which is closely related to that of Barge, Diamond, Hunton and Sadun \cite{BDHS10}. The descriptions of the PE and ePE homology groups that appear in this paper for the Penrose kite and dart tilings are made possible through this method of calculation. The method is directly applicable to a broad range of tilings, including `mixed substitution tilings' (see \cite{GahMal13}) but also non-Euclidean examples, such as the pentagonal tilings of Bowers and Stephenson mentioned above. The `approximant homology groups' of the computation and the `connecting maps' between them have a direct description in terms of the combinatorics of the star patches, making it highly amenable to computer implementation. In \cite{Gon11}, Gon\c{c}alves used the duals of these approximant chain complexes for a computation of the $K$-theory of the $C^*$-algebra of the stable equivalence relation of a substitution tiling. Our method of computation of the PE homology groups seems to confirm the observation there of a certain duality between these $K$-groups and the $K$-theory of the tiling space.

\subsection*{Organisation of Paper} In Section \ref{sec: Tilings and Tiling Spaces} we shall recall how one may associate to a Euclidean tiling $\mf{T}$ its translational hull $\Omega_\mf{T}^1$. When $\mf{T}$ has FLC, we also describe the presentation of $\Omega_\mf{T}^1$ as an inverse limit of approximants. In Section \ref{sec: Translational Pattern-Equivariance} we recall the PE cohomology of an FLC tiling $\mf{T}$, and how it may be identified with the \v{C}ech cohomology $\check{H}^\bullet(\Omega_\mf{T}^1)$ of the tiling space. We then introduce the PE homology of an FLC tiling and establish PE Poincar\'{e} duality between the PE cohomology and PE homology.

In Section \ref{sec: PE Homology and Rotations} we consider PE homology in the context of rotational symmetry. The ePE (co)homology groups are defined in Subsection \ref{subsec: Euclidean Pattern-Equivariance}, where we show, in Theorem \ref{thm: ePE PD}, that the ePE cohomology and ePE homology are Poincar\'{e} dual when taken over suitably divisible coefficients. In Subsection \ref{subsec: Restoring PD} we show how Poincar\'{e} duality for the ePE homology of a two-dimensional tiling is restored for $\Z$ coefficients by a simple modification of the ePE homology. The action of rotation on the PE cohomology of an FLC tiling, and its interaction with the ePE homology, is considered in Subsection \ref{subsec: Rotation Actions on Translational PE Cohomology}. In Subsection \ref{subsec: Generalising} we demonstrate how the techniques of Section \ref{sec: PE Homology and Rotations} may be naturally extended to a more general framework.

In Section \ref{sec: PE Homology of Hierarchical Tilings} we develop a method of computation of the PE homology for polytopal substitution tilings, close in spirit to the BDHS approach \cite{BDHS10}. In Subsection \ref{subsec: Generalisations} we explain how the method is modified to compute the ePE homology, and how it may be applied to more general settings, such as mixed substitution systems or to non-Euclidean examples.

\subsection*{Acknowledgements}

The author thanks John Hunton, Alex Clark, Lorenzo Sadun and Dan Rust for numerous helpful discussions, and the anonymous referee for their valuable suggestions.


\section{Tilings and Tiling Spaces} \label{sec: Tilings and Tiling Spaces}

\subsection{Cellular, Polytopal and Dual Tilings} \label{subsec: Cellular, Polytopal and Dual Tilings}

Recall that a CW complex is called \emph{regular} if the attaching maps of its cells may be taken to be homeomorphisms. A \emph{cellular tiling} of $\R^d$ shall be defined to be a pair $\mf{T} = (\mc{T},l)$ of a regular CW decomposition $\mc{T}$ of $\R^d$ along with a \emph{labelling} $l$ of $\mc{T}$, by which we mean a map from the cells of $\mc{T}$ to some set of `labels' $L$. We shall take \emph{cell} to mean a closed cell.  If the cells are convex polytopes then we call $\mf{T}$ a \emph{polytopal tiling}. For brevity, we will often refer to a cellular tiling as simply a \emph{tiling}, and a $d$-cell of a tiling as a \emph{tile}. A \emph{patch} of $\mf{T}$ is a finite subcomplex $\mc{P}$ of $\mc{T}$ together with the labelling restricted to $\mc{P}$. For a bounded set $U \subseteq \R^d$, we let $\mf{T}[U]$ be the patch supported on the set of tiles $t$ for which $t \cap U \neq \emptyset$.

Homeomorphisms of $\R^d$ act on tilings and patches in the obvious way. Two patches are called \emph{translation equivalent} if one is a translate of the other. The \emph{diameter} of a patch is defined to be the diameter of the support of its tiles. A tiling or, more generally, a collection of tilings, is said to have (\emph{translational}) \emph{finite local complexity} (\emph{FLC}) if, for any $r>0$, there are only finitely many patches of diameter at most $r$ up to translation equivalence. It is not difficult to see that a cellular tiling has FLC if and only if there are only finitely many translation classes of cells and the labelling function takes on only finitely many distinct values.

One may wish to consider other forms of decoration of $\R^d$, such as Delone sets, or tilings with overlapping or fractal tiles, and many of the concepts that we describe here have obvious analogues for them. However, when such a pattern has FLC it is always essentially equivalent to a polytopal tiling. In more detail, we say that $\mf{T}'$ \emph{is locally derivable from} $\mf{T}$ if there exists some $r>0$ for which, whenever $(S_1 \mf{T})[B_r] = (S_2 \mf{T})[B_r]$, then $(S_1 \mf{T}')[B_0] = (S_2\mf{T}')[B_0]$, where the $S_i$ are translations $x \mapsto x + t_i$ and $B_r$ is the closed ball of radius $r$ centred at the origin. The tilings $\mf{T}$ and $\mf{T}'$ are said to be \emph{mutually locally derivable} (\emph{MLD}) if each is locally derivable from the other. Loosely, this means that $\mf{T}$ and $\mf{T}'$ only differ in a very cosmetic sense, via locally defined redecoration rules. This concept was introduced in \cite{BSJ91}, along with the finer relation of \emph{S-MLD} equivalence, which takes into account general Euclidean isometries rather than just translations by replacing the translations $S_i$ in the definition of a local derivation above with Euclidean motions. FLC patterns (or even eFLC ones, see Subsection \ref{subsec: Rotational Tiling Spaces}) are always S-MLD to polytopal tilings via a Voronoi construction.

In the following sections we shall usually take our tilings to be polytopal. Since the properties of tilings of interest to us are preserved under S-MLD equivalence, this is not a harsh restriction. The major motivation for this choice is that some useful constructions may be defined for a polytopal complex, namely the barycentric subdivision and dual complex. In fact, it is sufficient for these constructions to use regular CW complexes as our starting point, but the most efficient way of dealing with this more general case is to pass to a combinatorial setting, which we do not cover in full detail here although shall outline in Subsection \ref{subsec: Generalising}.

For a polytopal tiling $\mf{T} = (\mc{T},l)$ we may construct the barycentric subdivision $\mc{T}_\Delta$ of its underlying CW decomposition geometrically, as follows. For each cell $c \in \mc{T}$, define $b(c)$, its \emph{barycentre}, to be the centre of mass of $c$ in its supporting hyperplane. We write $c_1 \preceq c_2$ for closed cells $c_1, c_2$ of $\mc{T}$ to mean that $c_1 \subseteq c_2$, and $c_1 \prec c_2$ if this inclusion is strict. The $k$-skeleton $\mc{T}_\Delta^k$ for $k = 0, \ldots, d$ is defined by taking as $k$-cells simplices which are the convex hulls of the vertices $\{b(c_0), b(c_1), \ldots, b(c_k)\}$ for a chain of cells $c_0 \prec c_1 \prec \cdots \prec c_k$ of $\mc{T}$ of length $k+1$. Such a cell may be labelled by the sequence of labels $(l(c_0),l(c_1),\ldots,l(c_k))$; we define the tiling $\mf{T}_\Delta = (\mc{T}_\Delta,l_\Delta)$, where $l_\Delta$ is the labelling of the cells of $\mc{T}_\Delta$ defined in this way. Assuming (without loss of generality, up to S-MLD equivalence) that cells of different dimension have different labels, it is easily verified that $\mf{T}$ and $\mf{T}_\Delta$ are S-MLD.

We may reconstruct $\mc{T}$ from its barycentric subdivision $\mc{T}_\Delta$ by identifying an open $k$-cell $c$ of $\mc{T}$ with the conglomeration of open simplicial cells corresponding to chains $c_0 \prec \cdots \prec c_j \prec c$ terminating in $c$. Flipping this process on its head, we obtain the \emph{dual complex} $\hat{\mc{T}}$. That is, we define an open $k$-cell of $\hat{\mc{T}}$ as the union of open simplicial cells corresponding to chains $c \prec c_0 \prec \cdots \prec c_j$ emanating from a $(d-k)$-cell $c$. Similarly to $\mc{T}_\Delta$, we may easily label $\hat{\mc{T}}$ so as to define a \emph{dual tiling} of $\mf{T}$ which is S-MLD to $\mf{T}_\Delta$, and hence also S-MLD to the original tiling $\mf{T}$. The dual tiling $\hat{\mf{T}}$ typically won't have convex polytopal cells, but it is cellular, owing to the piecewise linearity of the polytopal decomposition $\mc{T}$. The $k$-cells of $\mf{T}$ are naturally in bijection with the $(d-k)$-cells of $\hat{\mf{T}}$, and we have that $a \preceq b$ for cells of $\mf{T}$ if and only if $\hat{a} \succeq \hat{b}$ for the corresponding dual cells of $\hat{\mf{T}}$. A similar process would have worked for $\mc{T}$ only regular cellular. However, the decomposition of $\R^d$ defined by $\hat{\mc{T}}$ need not be cellular even for non-piecewise linear \emph{simplicial} complexes $\mc{T}$. Even so, the resulting dual decomposition $\hat{\mc{T}}$ still retains the analogous homological properties to a CW complex needed to define cellular homology (see \cite[Ch.\ 8.64]{Mun84}) and so the constructions and arguments to follow can, with little extra effort, be extended to this case.

\subsection{Tiling spaces}

To a tiling $\mf{T}$ one may associate a moduli space $\Omega_\mf{T}^1$ which, as a set, consists of tilings `locally indistinguishable' from $\mf{T}$. Let $S$ be a set of tilings of $\R^d$. We wish to endow $S$ with a geometry which expresses the intuitive idea that two tilings are close if, up to a small perturbation, the two agree to a large radius about the origin.

An approach which neatly applies to a large class of tilings, and captures this idea most directly, proceeds as follows. Let $H(\R^d)$ be the space of homeomorphisms of $\R^d$; these shall serve as our perturbations. Equipped with the compact--open topology, we may consider a neighbourhood $V \subseteq H(\R^d)$ of the identity as `small' if its elements only perturb points within a large distance from the origin of $\R^d$ a small amount. In this case, for $f \in V$ it is intuitive to think of $f \mf{T}$ as a small perturbation of $\mf{T}$. In fact, $\mf{T}$ should still be `close' to any other tiling $\mf{T}'$ so long as $\mf{T}[K]$ and $(f \mf{T}')[K]$ agree for some $K \subseteq \R^d$ containing a large neighbourhood of the origin. For a neighbourhood $V \subseteq H(\R^d)$ of $\operatorname{id}_{\R^d}$ and bounded $K \subseteq \R^d$, we define
\[
U(K,V) := \{(\mf{T}_1,\mf{T}_2) \in S \times S \mid \mf{T}_1[K] = (f \mf{T}_2)[K]\}.
\]
It is not difficult to verify that the resulting collection of sets $U(K,V)$ is a base for a uniformity on $S$. If the reader is unfamiliar with uniformities, the only important point here is that we have a uniform notion of tilings being `close': \emph{$\mf{T}_1$ is considered as `close' to $\mf{T}_2$, as judged by $K$ and $V$, whenever $(\mf{T}_1,\mf{T}_2) \in U(K,V)$}. With $K$ a large neighbourhood of the origin and $V$ a set of homeomorphisms moving points only a very small amount in the vicinity of $K$, we recover our intuitive notion of $\mf{T}_1$ and $\mf{T}_2$ being close, when they agree to a large radius up to a small perturbation. The above construction easily generalises to other decorations of $\R^d$, such as Delone sets, and also tilings with infinite label sets which are equipped with a metric (see \cite{PrFSad14ILC}).

For a tiling $\mf{T}$ we define the \emph{translational hull} or \emph{tiling space} as
\[
\Omega^1_\mf{T} := \overline{\{\mf{T} + x \mid x \in \R^d\}},
\]
where the completion is taken with respect to the uniformity defined above. In the case that $\mf{T}$ has FLC, two patches agree up to a small perturbation when they agree up to a small translation. So in this case the sets
\[
U(K,\epsilon) \coloneqq \{(\mf{T}_1,\mf{T}_2) \mid \mf{T}_1[K] = (\mf{T}_2 + x)[K] \text{ for } \|x\| \leq \epsilon\}
\]
serve as a base for our uniformity, where the $K \subseteq \R^d$ are bounded and $\epsilon > 0$. Loosely, two tilings are `close' if and only if their central patches agree to a large radius (parametrised by $K$) up to a small translation (parametrised by $\epsilon$). It is not difficult to show that the tiling space $\Omega^1_\mf{T}$ is a compact space whose points may be identified with those tilings whose patches are translates of patches of $\mf{T}$. So one may take as basic open neighbourhoods of a tiling $\mf{T}_1$ in $\Omega^1_\mf{T}$ the \emph{cylinder sets}
\[
C(R,\epsilon,\mf{T}_1) \coloneqq \{\mf{T}_2 \in \Omega^1_\mf{T} \mid \mf{T}_1[B_R] = (\mf{T}_2 + x)[B_R] \text{ for } \|x\| \leq \epsilon\}
\]
of tilings of the hull which, up to a translation of at most $\epsilon$, agree with $\mf{T}_1$ to radius $R$.

\subsection{Inverse Limit Presentations} \label{sec: Inverse Limit Presentations}

Another simplification granted to us by finite local complexity is that the tiling space $\Omega^1_\mf{T}$ may be presented as an inverse limit of CW complexes $\Gamma^1_i$, following G\"{a}hler's (unpublished) construction; see \cite{Sad08} for details, and also the alternative approach of Barge, Diamond, Hunton and Sadun \cite{BDHS10}. Inductively define the \emph{$i$-corona} of a tile as follows: the $0$-corona of a tile $t$ is the patch whose single tile is $t$; for $i \in \mathbb{N}$, the $i$-corona is the patch of tiles which have non-empty intersection with the $(i-1)$-corona. That is, one constructs the $i$-corona of $t$ by taking $t$ and then iteratively appending neighbouring tiles $i$ times. For $x,y \in \R^d$, write $x \sim_i y$ to mean that there are two tiles $t_x$ and $t_y$ of $\mf{T}$ containing $x$ and $y$, respectively, for which the $i$-corona of $t_x$ is equal to the $i$-corona of $t_y$, up to a translation taking $x$ to $y$. This is typically not an equivalence relation, so we define the \emph{approximant} $\Gamma^1_i$ to be the quotient of $\R^d$ by the transitive closure of the relation $\sim_i$. More intuitively, we form $\Gamma^1_i$ by taking a copy of the central tile from each translation class of $i$-corona, glueing them along their boundaries according to how they can meet in the tiling. We define the \emph{$i$-corona} of a lower dimensional cell to be the intersection of $i$-coronas of the tiles which it is contained in. An alternative way of defining approximants, which avoids taking a transitive closure (although identifies more points of $\R^d$ at each level), is to identify cells of the tiling which share the same $i$-coronas, up to a translation. Each approximant naturally inherits a cellular decomposition from that of the tiling.

For $i \leq j$, cells of $\mf{T}$ identified in $\Gamma^1_j$ are also identified in $\Gamma^1_i$, so we have `forgetful' cellular quotient maps $\pi_{i,j} \co \Gamma^1_j \rightarrow \Gamma^1_i$. The inverse limit of this projective system
\[
\varprojlim (\Gamma^1_i, \pi_{i,j}) \coloneqq \{(x_i)_{i \in \mathbb{N}_0} \in \prod_{i=0}^\infty \Gamma^1_i \mid \pi_{i,j}(x_j) = x_i\}
\]
is homeomorphic to the tiling space $\Omega^1_\mf{T}$. The central idea here is that a point of $\Gamma^1_i$ describes how to tile a neighbourhood of the origin, where the sizes of these neighbourhoods increase with $i$. An element of the inverse limit space then corresponds to a consistent sequence of choices of larger and larger patches about the origin, so it defines a tiling. Any two tilings which are `close' correspond to points of the inverse limit which are `close' on an approximant $\Gamma^1_i$ for large $i$, and vice versa.

\section{Translational Pattern-Equivariance} \label{sec: Translational Pattern-Equivariance}

\subsection{Identifying \v{C}ech with PE Cohomology} Locally, the tiling space of an FLC tiling has a product structure of cylinder sets $U \times C$ where $U$ is an open subset of $\R^d$, corresponding to small translations, and $C$ is a totally disconnected space, corresponding to different  ways of completing a finite patch to a full tiling. Globally, $\Omega^1_\mf{T}$ is a torus bundle with totally disconnected fibre \cite{SadWil03}. Many classical invariants---homotopy groups and singular (co)homology groups, for example---are ill-suited to studying $\Omega^1_\mf{T}$ when $\mf{T}$ is non-periodic, in which case this space is not locally connected. A commonly employed topological invariant with which to study tiling spaces is \v{C}ech cohomology $\check{H}^\bullet(-)$. We shall not cover its definition here (see \cite[Chapter\ 2.10]{BottTu82}), although we recall two important features of it:
\begin{enumerate}
	\item \v{C}ech cohomology is naturally isomorphic to singular cohomology on the category of spaces homotopy equivalent to CW complexes and continuous maps.
	\item For a projective system $(\Gamma_i,\pi_{i,j})$ of compact, Hausdorff spaces $\Gamma_i$, we have an isomorphism $\check{H}^\bullet(\varprojlim (\Gamma_i,\pi_{i,j})) \cong \varinjlim (\check{H}^\bullet(\Gamma_i),\pi_{i,j}^*)$.
\end{enumerate}

Pattern-equivariant cohomology is a tool designed to give intuitive descriptions of the \v{C}ech cohomology of tiling spaces. It was first defined by Kellendonk and Putnam in \cite{KelPut06} (see also \cite{Kel03}) where they showed that it is isomorphic to the \v{C}ech cohomology of the tiling space taken over $\R$ coefficients. It is constructed by restricting the de Rham cochain complex of $\R^d$ of smooth forms to a sub-cochain complex of forms which, loosely, are determined pointwise by the local decoration of the underlying tiling to some bounded radius.

A second approach, introduced by Sadun in \cite{Sad07}, is to use cellular cochains, and has the advantage of generalising to arbitrary Abelian coefficients. Let $\mf{T} = (\mc{T},l) $ be a cellular tiling (recall that $\mc{T}$ is the underlying cell complex of $\mf{T}$). Denote by $C^\bullet(\mc{T})$ the cellular cochain complex of $\mc{T}$;
\[
C^\bullet(\mc{T}) : = 0 \rightarrow C^0(\mc{T}) \xrightarrow{\delta^0} C^1(\mc{T}) \xrightarrow{\delta^1} \cdots \xrightarrow{\delta^{d-1}} C^d(\mc{T}) \rightarrow 0,
\]
where each $C^k(\mc{T})$ is the group of cellular $k$-cochains and $\delta^k$ is the degree $k$ cellular coboundary map. A cellular $k$-cochain $\psi$ is a function which assigns to each orientation $\omega_c$ of $k$-cell $c$ an integer, satisfying $\psi(\omega_c^+) = - \psi(\omega_c^-)$ for opposite orientations $\omega_c^+$ and $\omega_c^-$ of a cell $c$. Of course, choosing an orientation for each $k$-cell induces an isomorphism $C^k(\mc{T}) \cong \prod _{k \text{-cells}} \Z$. Choose orientations for the $k$-cells so that $\omega_c + x = \omega_{c+x}$ whenever $c$ and $c+x$ are cells of $\mc{T}$, where $\omega_c + x$ is the orientation on $c+x$ induced from $\omega_c$ by translation. Write $\psi(c) \coloneqq \psi(\omega_c)$, where $\omega_c$ is the chosen orientation of $c$. A cochain $\psi$ is called \emph{pattern-equivariant} (\emph{PE}) if there exists some $i \in \mathbb{N}_0$ for which $\psi(c_1) = \psi(c_2)$ whenever $c_1$ and $c_2$ have identical $i$-coronas, up to a translation taking $c_1$ to $c_2$. 

It is easy to check that the coboundary of a PE cochain is PE. Define $C^\bullet(\mf{T}^1)$ to be the sub-cochain complex of $C^\bullet(\mc{T})$ consisting of PE cochains. Its cohomology $H^\bullet(\mf{T}^1)$ is called the \emph{pattern-equivariant cohomology of} $\mf{T}$.

A cellular cochain $\psi \in C^k(\mc{T})$ is PE if and only if it is a pullback cochain from some approximant, that is, if $\psi = \pi_i^*(\widetilde{\psi})$ where $\widetilde{\psi} \in C^k(\Gamma^1_i)$ is a cellular cochain on an approximant and $\pi_i$ is the (cellular) quotient map $\pi_i \co \R^d \rightarrow \Gamma^1_i$ defining $\Gamma^1_i$. This fact, in combination with the description of the tiling space $\Omega^1_\mf{T}$ as an inverse limit of G\"{a}hler complexes and the two features of \v{C}ech cohomology given above, leads to the proof in \cite{Sad07} of the following:

\begin{theorem}[\cite{Sad07}] \label{thm: Cech = PE} The PE cohomology $H^\bullet(\mf{T}^1)$ of an FLC tiling $\mf{T}$ is isomorphic to the \v{C}ech cohomology $\check{H}^\bullet(\Omega^1_\mf{T})$ of its tiling space. \end{theorem}

\subsection{PE Homology and Poincar\'{e} Duality}

We shall now define the PE \emph{homology} groups of a cellular tiling $\mf{T}$. The construction runs almost identically to the construction of the cellular PE cohomology groups above, but where we took cellular coboundary maps before we shall take instead cellular \emph{boundary} maps. In more detail, let $C_\bullet^{\text{BM}}(\mc{T})$ denote the \emph{cellular Borel--Moore chain complex}
\[
C_\bullet^{\text{BM}}(\mc{T}) \coloneqq 0 \leftarrow C_0^{\text{BM}}(\mc{T}) \xleftarrow{\partial_1} C^{\text{BM}}_1(\mc{T}) \xleftarrow{\partial_2} \cdots \xleftarrow{\partial_d} C^{\text{BM}}_d(\mc{T}) \leftarrow 0.
\]
The chain groups $C_k^{\text{BM}}(\mc{T})$ are canonically isomorphic to the cochain groups $C^k(\mc{T})$. That is, up to a choice of orientations for the $k$-cells, a cellular Borel--Moore chain $\sigma \in C_k^{\text{BM}}(\mc{T})$ is given by a choice of integer for each $k$-cell. But we think of its elements as possibly infinite, or non-compactly supported cellular \emph{chains}. The boundary maps $\partial_k$ are the linear extension to these chain groups of the cellular boundary maps of the standard cellular chain complex of $\mc{T}$.

Pattern-equivariance of a chain $\sigma \in C_k^{\text{BM}}(\mc{T})$ is defined identically to that of a cochain. That is, $\sigma$ is PE if there exists some $i \in \mathbb{N}_0$ for which, for any two $k$-cells $c_1, c_2$ of $\mc{T}$ with identical $i$-coronas in $\mf{T}$ up to a translation, $c_1$ and $c_2$ have the same coefficient in $\sigma$. It is easy to see that if $\sigma$ is PE then $\partial(\sigma)$ is also PE. Restricting to PE cellular Borel--Moore chains we obtain a sub-chain complex $C_\bullet(\mf{T}^1)$ of $C^{\text{BM}}_\bullet(\mc{T})$ whose homology $H_\bullet(\mf{T}^1)$ we shall call the \emph{pattern-equivariant homology of} $\mf{T}$. So the elements of the PE homology are represented by, typically, non-compactly supported cellular \emph{cycles} (chains with trivial boundary) where two PE cycles $\sigma_1, \sigma_2$ are homologous if $\sigma_1 = \sigma_2 + \partial(\tau)$ for some PE chain $\tau$.

These homology groups certainly have a highly geometric definition, but what do they measure? Through a Poincar\'{e} duality argument, we may in fact identify them with the (re-indexed) PE cohomology groups and thus, in light of Theorem \ref{thm: Cech = PE}, with the \v{C}ech cohomology groups of the tiling space:

\begin{theorem} \label{thm: PE PD} For a polytopal tiling $\mf{T}$ of $\R^d$ of finite local complexity, we have PE Poincar\'{e} duality $H^\bullet(\mf{T}^1) \cong H_{d - \bullet}(\mf{T}^1)$.\end{theorem}

\begin{proof} Classical Poincar\'{e} duality provides an isomorphism of complexes
\[
- \cap \Gamma \co C^\bullet(\mc{T}) \xrightarrow{\cong} C^{\text{BM}}_{d-\bullet}(\hat{\mc{T}}),
\]
via the cap product with a cellular Borel--Moore fundamental class $\Gamma$, a $d$-cycle for which each oriented $d$-cell has coefficient either $+1$ or $-1$, pairing orientations of $k$-cells with orientations of their dual $(d-k)$-cells. Here, $\mc{T}$ is the underlying cell complex of $\mf{T}$ and $\hat{\mc{T}}$ is its dual complex. By definition, a (co)chain is PE whenever it assigns coefficients to cells in a way which only depends locally on the tiling to some bounded neighbourhood of that cell. The fundamental class $\Gamma$ is also PE. Since the cap product (and here, its inverse) is defined in a local manner, and $\mf{T}$ and $\hat{\mf{T}}$ are MLD, a $k$-cochain $\psi$ of $\mf{T}$ is PE if and only if its dual $(d-k$)-chain $\psi \cap \Gamma$ of $\hat{\mf{T}}$ is PE. So $- \cap \Gamma$ restricts to an isomorphism between PE complexes
\[
- \cap \Gamma \co C^\bullet(\mf{T}^1) \xrightarrow{\cong} C_{d-\bullet}(\hat{\mf{T}^1}).
\]

The barycentric tiling $\mf{T}_\Delta$ refines both $\mf{T}$ and $\hat{\mf{T}}$, and is MLD to both. As one may expect, taking such a refinement does not effect PE (co)homology. This shall be shown, in a more general setting, in Lemma \ref{lem: refinement}. Precisely, we have quasi-isomorphisms $\iota \co C_\bullet(\mf{T}^1) \rightarrow C_\bullet(\mf{T}_\Delta^1)$ and $\hat{\iota} \co C_\bullet(\hat{\mf{T}}^1) \rightarrow C_\bullet(\mf{T}_\Delta^1)$; recall that a quasi-isomorphism is a (co)chain map which induces an isomorphism on (co)homology. In summation we have the following diagram of quasi-isomorphisms
\[
C^\bullet(\mf{T}^1) \xrightarrow{- \cap \Gamma} C_{d-\bullet}(\hat{\mf{T}}^1) \xrightarrow{\hat{\iota}} C_{d-\bullet}(\mf{T}_\Delta^1) \xleftarrow{\iota} C_{d-\bullet}(\mf{T}^1)
\]
from which PE Poincar\'{e} duality $H^\bullet(\mf{T}^1) \cong H_{d - \bullet}(\mf{T}^1)$ follows.
\end{proof}

\begin{figure}
\begin{center}
\includegraphics[scale=0.6]{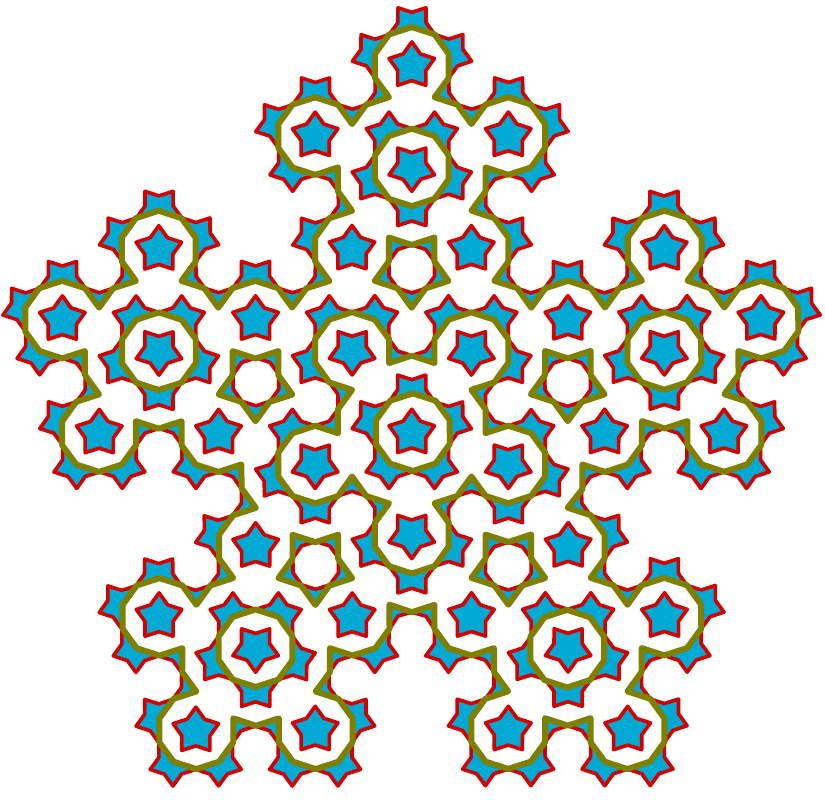} \caption{The PE $1$-cycle $\rho'$ of a Penrose kite and dart tiling (red), with the analogous chain of the supertiling (green) and a PE $2$-chain (blue) whose boundary relates the two.}
\label{fig: Penrose}
\end{center}
\end{figure}

\begin{exmp} \label{ex: Penrose PE} Let $\mf{T}$ be a Penrose kite and dart tiling. The \v{C}ech cohomology of the translational hull of the Penrose tilings was first calculated in \cite{AndPut98} (although see also the earlier closely related $K$-theoretic calculations of Kellendonk in the groupoid setting \cite{Kel97}). In Section \ref{sec: PE Homology of Hierarchical Tilings} we provide a different way of computing these groups which, as a direct by-product, provides us with explicit descriptions of the generators in terms of PE chains. Consistently with previous calculations, we find that
\[
H_{2-k}(\mf{T}^1) \cong H^k(\mf{T}^1) \cong \check{H}^k(\Omega^1_\mf{T}) \cong \begin{cases}
	\Z \text{ for } k=0;\\
	\Z^5 \text{ for } k=1;\\
	\Z^8 \text{ for } k=2.\\
	\end{cases}
\]
Let $P_c$ be a pair of a patch $P$ from $\mf{T}$ along with a choice of oriented $k$-cell $c$ from $P$, taken up to translation equivalence. We have an associated PE \emph{indicator $k$-chain} $\cf(P_c) \in C_k(\mf{T}^1)$ for which each $k$-cell of $\mf{T}$ has coefficient one when it is contained in an ambient patch for which the pair agrees with $P_c$ up to translation, and all other cells have coefficient zero. The degree zero PE homology group $H_0(\mf{T}^1)$ for a Penrose kite and dart tiling may be freely generated by indicator $0$-chains $\cf(P_v) \in C_0(\mf{T}^1)$, where each $P_v$ is one of the vertex-stars of $\mf{T}$, paired with its central vertex. The full list of possible translation classes of such patches, up to rotation by some $2 \pi k /10$, are given (and named, according to Conway's notation) in Figure \ref{fig: VaETypes}. As an example of a choice of elements freely generating $H_0(\mf{T}^1)$, we may choose two queen vertices, one oriented as in Figure \ref{fig: VaETypes} and the other a $2 \pi/10$ rotate of it, and six king vertices, each a $2\pi k/10$ rotate of that of Figure \ref{fig: VaETypes} with $k=0,\ldots,5$.

We shall go into more detail on generators for $H_0(\mf{T}^1)$ in Example \ref{ex: Penrose results}. In degree one, there are two particularly beautiful cycles that we wish to discuss here. There is a PE $1$-cycle $\rho'$ given by running along the bottoms of the dart tilings, with $1$-cells oriented to point to, say, the right when the darts are oriented to point upwards. The resulting cycle is illustrated in red in Figure \ref{fig: Penrose}, where we have removed cell orientations and the $1$-skeleton of the tiling to decrease clutter. The extra embellishments of the figure shall be discussed in Section \ref{sec: PE Homology of Hierarchical Tilings}; there is a green $1$-cycle for the analogous chain of the supertiling of $\mf{T}$, along with a blue PE $2$-chain whose boundary relates the two. As one can see, $\rho'$ forms a disjoint union of clockwise and anticlockwise running loops. Interestingly, deducing which of these two options is the case at some cell of a loop requires consideration of arbitrarily large patches; in fact, for specific kite and dart tilings there exists a single infinite fractal-like path along dart tiles. But $\rho'$ is not a generator, there exists another PE $1$-cycle $\rho$ for which $[\rho'] = 2[\rho]$ in $H_1(\mf{T}^1)$. The loops of Figure \ref{fig: Penrose} come in two types: ones where the darts along the loops are $2\pi k /10$ rotates of an upwards pointing dart tile with $k$ odd, and ones where the darts are even rotates. The $1$-cycle $\rho$ is given by restricting $\rho'$ to those loops in one of these two parities; both choices are homologous and equal to $[\rho']/2$ in $H_1(\mf{T}^1)$.

A second generator $\nu_0$ for $H_1(\mf{T}^1)$ is depicted in Figure \ref{fig: Ammann Chain}. The cycle arranges as a union of infinite paths along the $1$-skeleton which closely approximate the \emph{Ammann bars} \cite[Chapter 10.6]{GruShe87} of the supertiling of $\mf{T}$, illustrated in the figure in green. There are ten further chains $\nu_k$ defined by $2\pi k /10$ rotates of $\nu_0$ (see Subsection \ref{subsec: Rotation Actions on Translational PE Cohomology}). We calculate that $H_1(\mf{T}^1)$ is freely generated by the homology classes of $\rho$, $\nu_0$, $\nu_1$, $\nu_2$ and $\nu_3$; every other PE $1$-cycle is equal, up to a PE $2$-boundary, to a linear combination of these cycles. It turns out that $\nu_4 \simeq -\nu_0+\nu_1-\nu_2+\nu_3$ in homology. This formula is unsurprising following the observation that one may associate each $\nu_k$ with a direction given by a tenth root of unity, and we have the identity $\sum_{k=0}^4 (-1)^k \exp(2 \pi i \cdot k /10) = 0$. \end{exmp}


\section{PE Homology and Rotations} \label{sec: PE Homology and Rotations}

In the previous section we showed that topological invariants of tiling spaces may be described in a highly geometric way, using infinite cellular chains on the tiling. However, PE homology is essentially just offering a different perspective on the generators of the PE cohomology here. As we shall see in Section \ref{sec: PE Homology of Hierarchical Tilings}, PE homology does provide a valuable alternative insight into the calculation of these invariants for hierarchical tilings. In this section, we shall show that PE homology provides a new invariant to the PE cohomology when one considers these invariants in the context of rotational symmetries.

\subsection{Rotational Tiling Spaces} \label{subsec: Rotational Tiling Spaces} Let $\text{S}E(d) \cong \R^d \rtimes \text{SO}(d)$ denote the transformation group of orientation-preserving isometries of $\R^d$, elements of $\text{S}E(d)$ shall be called \emph{rigid motions}. There are two topological spaces naturally associated to a tiling $\mf{T}$ of $\R^d$ which incorporate the action of $\text{S}E(d)$ on $\mf{T}$. The first, defined analogously to the translational hull $\Omega_\mf{T}^1$, is the \emph{Euclidean hull}
\[
\Omega_\mf{T}^{\text{rot}} \coloneqq \overline{ \{f \mf{T} \mid f \in \text{S}E(d)\}}.
\]
It follows directly from the definitions that the special Euclidean group $\text{S}E(d)$ acts uniformly on the Euclidean orbit of $\mf{T}$, and so this action extends to the entire Euclidean hull. In particular, the subgroup $\text{SO}(d) \leq \text{S}E(d)$ of rotations at the origin acts on $\Omega_\mf{T}^{\text{rot}}$. The second space, the one which we shall concentrate on in this section, is the quotient
\[
\Omega_\mf{T}^0 \coloneqq \Omega_\mf{T}^{\text{rot}} / \text{SO}(d).
\]

We shall say that $\mf{T}$ has \emph{Euclidean finite local complexity} (\emph{eFLC}, for short) if, for every $r>0$, there exist only finitely many patches of diameter at most $r$ up to rigid motion. Interesting examples of tilings which have eFLC, but not translational FLC, are the Conway--Radin pinwheel tilings of $\R^2$, whose tiles are all rigid motions of a $(1,2,\sqrt{5})$ triangle, or its reflection, but are found in the tiling pointing in infinitely many directions. Much of what can be said for FLC tilings and their translational hulls has an analogue for eFLC tilings and their Euclidean hulls. In particular, for an eFLC tiling $\mf{T}$, its Euclidean hull $\Omega_\mf{T}^{\text{rot}}$ is a compact space whose points may be identified with those tilings whose patches are rigid motions of the patches of $\mf{T}$. The space $\Omega_\mf{T}^0$ is then the quotient of $\Omega_\mf{T}^{\text{rot}}$ given by identifying tilings which differ by a rotation at the origin. One may define inverse limit presentations of these spaces in a similar way to the construction of the G\"{a}hler complexes, which is tantamount to being able to define pattern-equivariant cohomology.

To explain this further, we now focus on the space $\Omega_\mf{T}^0$. For $i \in \mathbb{N}_0$ we define CW complexes $\Gamma_i^0$ analogously to the complexes of the translational setting, replacing translations with rigid motions. For example, we may define the complexes $\Gamma_i^0$ by identifying cells $c_1$, $c_2$ of $\mf{T}$ via rigid motions which take $c_1$ to $c_2$, and the $i$-corona of $c_1$ to the $i$-corona of $c_2$. The CW complexes $\Gamma_i^0$, along with the `forgetful maps' between them, define a projective system whose inverse limit is homeomorphic to $\Omega_\mf{T}^0$.

It may be the case that cells of $\mf{T}$ have \emph{non-trivial isotropy}, that is, there may be cells $c$ whose $i$-coronas are preserved by rigid motions mapping $c$ to itself non-trivially, which will cause points of $c$ to be identified in the quotient spaces $\Gamma_i^0$. Given a cell $c \in \mc{T}$, the rigid motions mapping $c$ to itself and preserving its $i$-corona is a group, which we call the \emph{isotropy}, and denote by $\text{Iso}(c,i)$. Write $\widetilde{\text{Iso}}(c,i)$, the \emph{cell isotropy}, to denote the group of transformations of $\text{Iso}(c,i)$ restricted to $c$.

The cell isotropy groups of the barycentric subdivision $\mf{T}_\Delta$ are always trivial. Indeed, a barycentric cell is determined by its vertices, which are determined by a chain of incidences $c_0 \prec \cdots \prec c_k$ in $\mc{T}$, and a rigid motion taking a barycentric cell to itself is determined by the map restricted to these vertices. A non-trivial map on such a vertex set would have to correspond to a rigid motion taking $b(c_i)$ to some $b(c_j)$ with $i \neq j$, which cannot be the case since $c_i \neq c_j$ have distinct dimensions, and thus, by assumption, distinct labels.

\subsection{Euclidean Pattern-Equivariance} \label{subsec: Euclidean Pattern-Equivariance}

A cellular cochain $\psi \in C^k(\mc{T})$ shall be called \emph{Euclidean pattern-equivariant} (\emph{ePE}) if there exists some $i \in \mathbb{N}_0$ for which $\psi(\omega_c) = \psi(f_*(\omega_c)) $ whenever $f$ is a rigid motion taking the $i$-corona of a $k$-cell $c$ to the $i$-corona of some other $k$-cell; here, $\omega_c$ is an orientation on $c$ and $f_*(\omega_c)$ is the push-forward of this orientation to the cell $f(c)$. In the case that the cells of $\mf{T}$ have trivial cell isotropy, one may consistently orient the cells of $\mf{T}$, and this definition then just says that there exists some $i$ for which $\psi$ is constant on cells which have identical $i$-coronas up to rigid motion. The coboundary of an ePE cochain is ePE, so we have a cochain complex $C^\bullet(\mf{T}^0)$ defined by restricting $C^\bullet(\mc{T})$ to ePE cochains. Taking the cohomology of this cochain complex, we define the \emph{ePE cohomology} $H^\bullet(\mf{T}^0)$.

One may follow the proof from \cite{Sad07} of Theorem \ref{thm: Cech = PE} almost word-for-word, replacing the G\"{a}hler complexes $\Gamma^1_i$ by the complexes $\Gamma_i^0$, to obtain the following:

\begin{theorem}\label{thm: Cech = ePE} Let $\mf{T}$ be an eFLC tiling whose cell isotropy groups $\widetilde{\text{Iso}}(c,i)$ are trivial for some $i \in \mathbb{N}_0$. Then the ePE cohomology $H^\bullet(\mf{T}^0)$ is isomorphic to the \v{C}ech cohomology $\check{H}^\bullet(\Omega^0_\mf{T})$.\end{theorem}

We define Euclidean pattern-equivariance for cellular Borel--Moore chains identically as for cochains. Restricting $C_\bullet^{\text{BM}}(\mc{T})$ to ePE chains we thus define the \emph{ePE chain complex} $C_\bullet(\mf{T}^0)$ and its homology, the \emph{ePE homology} $H_\bullet(\mf{T}^0)$.

The proof of PE Poincar\'{e} duality in Theorem \ref{thm: PE PD} essentially relied on two simple observations:
\begin{enumerate}
	\item The classical Poincar\'{e} duality isomorphism $- \cap \Gamma$, given by taking the cap product with a Borel--Moore fundamental class, restricts to a cochain isomorphism $- \cap \Gamma \co C^\bullet(\mf{T}^1) \rightarrow C_{d-\bullet}(\hat{\mf{T}}^1)$ between the PE cohomology of $\mf{T}$ and the PE homology of its dual tiling.
	\item The refinement of a tiling to its barycentric subdivision does not effect PE homology.
\end{enumerate}
Step $(1)$ will still hold for the ePE complexes, we have a Poincar\'{e} duality isomorphism $- \cap \Gamma \co C^\bullet (\mf{T}^0) \rightarrow C_{d-\bullet}(\hat{\mf{T}}^0)$ between ePE cochains of a tiling and ePE chains of its dual tiling. However, step $(2)$ will not generally hold when restricting to ePE (co)chains if our tiling has non-trivial cell isotropy. One would like to refuse taking the ePE (co)homology of a tiling which has cells of non-trivial isotropy by replacing it with the barycentric subdivision. Unfortunately, our hand is forced, since in the presence of non-trivial \emph{patch} isotropy one of $\mf{T}$ or $\hat{\mf{T}}$ will have non-trivial \emph{cell} isotropy.

The next lemma determines to what extent one may expect refinement to preserve ePE homology and cohomology. Thus far, our invariants have been taken over $\Z$ coefficients. For a unital ring $G$ we may consider the cochain complex $C^\bullet(\mc{T};G)$ of cellular cochains which assign to oriented cells elements of $G$, and similarly for $C_\bullet^{\text{BM}}(\mc{T};G)$. We may restrict these complexes to PE and ePE (co)chains, and denote the corresponding (co)homology by $H^\bullet(\mf{T}^1;G)$, etc. We say that $G$ has \emph{division by $n$} if $n \cdot 1_G$ has a multiplicative inverse in $G$, where $1_G$ is the multiplicative identity in $G$ and $n \cdot g \coloneqq \underbrace{g + g + \cdots + g}_{n \text{ times}}$.

\begin{lemma} \label{lem: refinement} Let $\mf{T}$ be a polytopal tiling with eFLC and $G$ be a unital ring. If for some $K \in \mathbb{N}_0$ the coefficient ring $G$ has division by $\# \widetilde{\text{Iso}}(c,K)$ for every cell $c$ of $\mf{T}$, then there exist quasi-isomorphisms
\[
\begin{array}{l}
\iota \co C_\bullet(\mf{T}^0;G) \rightarrow C_\bullet(\mf{T}_\Delta^0;G); \\
\kappa \co C^\bullet(\mf{T}_\Delta^0;G) \rightarrow C^\bullet(\mf{T}^0;G).
\end{array}
\]
The analogous statement holds for the refinement of the dual tiling $\hat{\mf{T}}$ to $\mf{T}_\Delta$.
\end{lemma}

\begin{proof} We shall prove the existence of the homology quasi-isomorphism, the proof for cohomology is similar. An \emph{elementary chain} is a chain which assigns coefficient $1$ to some oriented cell and $0$ to all other cells. We have a chain map
\[
\iota \co C^{\text{BM}}_\bullet(\mc{T};G) \rightarrow C^{\text{BM}}_\bullet(\mc{T}_\Delta;G)
\]
which is defined by sending an elementary chain to the corresponding barycentric chain with coefficient $1$ on barycentric $k$-cells contained in $c$, suitably oriented with respect to $c$, and $0$ on all other cells. It is easy to see that $\iota$ restricts to ePE chains and we claim that it is a quasi-isomorphism.

To show that $\iota$ is surjective on homology, let $\sigma \in C_k(\mf{T}^0_\Delta)$ be an ePE cycle of the barycentric subdivision; $\sigma$ is in the image of $\iota$ if and only if it is supported on the $k$-skeleton of $\mf{T}$. If $k = d$, then $\sigma$ is already supported on the $d$-skeleton, so suppose that $k < d$.

Whilst $\sigma$ need not be in the image of $\iota$ at the chain level, there exists some $\tau$ for which $\sigma + \partial(\tau)$ is. To construct $\tau$, we firstly find an ePE chain $\tau(d)$ for which $\sigma + \partial(\tau(d))$ is supported on the $(d-1)$-skeleton. Having the same $i$-corona up to rigid motion is an equivalence relation on the cells of $\mf{T}$ for every $i \in \mathbb{N}_0$. By Euclidean pattern-equivariance of $\sigma$, there exists some $i$ for which, whenever two $d$-cells $c_1, c_2$ of $\mf{T}$ have identical $i$-coronas up to a rigid motion $f$, then $f$ sends $\sigma$ restricted at $c_1$ to its restriction at $c_2$. We may choose $i \geq K$; note that, since $\widetilde{\text{Iso}}(c,i)$ is a subgroup of $\widetilde{\text{Iso}}(c,K)$, we have that $\# \widetilde{\text{Iso}}(c,i)$ divides $\# \widetilde{\text{Iso}}(c,K)$, so the coefficient group $G$ has division by $\#\widetilde{\text{Iso}}(c,i)$ for every cell $c$ of $\mf{T}$.

For each equivalence class of $d$-cell, choose a representative $c$ and a barycentric $(k+1)$-chain $\tau_c$ supported on $c$ for which $\sigma + \partial(\tau_c)$ is supported outside of the interior of $c$; by the homological properties of cells of a CW decomposition, we may find such a chain. Define the $(k+1)$-chain $\tau'_c$ by copying $\tau_c$ to every cell equivalent to $c$, via every rigid motion which preserves the $i$-corona of $c$. We define $\tau(d)  \coloneqq \sum \tau_c' / (\#\widetilde{\text{Iso}}(c,i))$, where the sum is taken over every equivalence class of $d$-cell.

The chain $\tau(d)$ is ePE by construction, and we claim that $\sigma + \partial(\tau(d))$ is supported on the $(d-1)$-skeleton. Indeed, let $c$ be a chosen representative of $d$-cell; we have that $\partial(\tau_c) = -\sigma_c$, where $\sigma_c$ is the restriction of $\sigma$ to the interior barycentric cells of $c$. By our assumption on $\sigma$ being ePE, for any $f \in \widetilde{\text{Iso}}(c,i)$ we have that $f_*(\partial(\tau_c)) =  f_*(-\sigma_c) = -\sigma_c$. Hence, the restriction of $\partial(\tau(d))$ to the interior of $c$ is given by:
\[
\sum_{f \in \widetilde{\text{Iso}}(c,i)} f_*(\partial(\tau_c))/(\# \widetilde{\text{Iso}}(c,i)) = \sum_{f \in \widetilde{\text{Iso}}(c,i)} -\sigma_c/(\# \widetilde{\text{Iso}}(c,i)) = -\sigma_c.
\]
By construction of $\tau(d)$, the same is true at every other $d$-cell equivalent to $c$, and by our assumption on $\sigma$ being ePE it follows that $\sigma + \partial(\tau(d))$ is supported on the $(d-1)$-skeleton.

We may continue this procedure down the skeleta. That is, we may construct in an analogous way ePE chains $\tau(d), \tau(d-1), \ldots, \tau(k+1)$ for which $\sigma + \partial (\sum_{m=n}^d \tau(m))$ is supported on the $n$-skeleton of $\mf{T}$. It follows that $\sigma + \partial (\sum_{m=k+1}^d \tau(m))$ is in the image of $\iota$, so $\iota_*$ is surjective on homology.

Showing injectivity of $\iota_*$ is analogous (indeed, the above is really just a relative homology argument applied to the filtration of the skeleta). Suppose that $\iota(\sigma) = \partial(\tau)$ for $\sigma \in C_k(\mf{T}^0)$ and $\tau \in C_{k+1}(\mf{T}^0_\Delta)$. Then $\tau$ is ePE and has boundary in the $k$-skeleton of $\mf{T}$. We may construct ePE $(k+2)$-chains $\nu(d),\nu(d-1),\ldots,\nu(k+2)$, analogously to above, for which $\tau + \partial(\sum_{m=k+2}^d \nu(m))$ is contained in the $(k+1)$-skeleton. So there is an ePE chain $\tau'$ with $\iota(\tau') = \tau + \partial(\sum_{m=k+2}^d \nu(m))$. It follows from $\partial (\iota(\tau')) = \partial(\tau) = \iota(\sigma)$ that $\sigma = \partial(\tau')$ represents zero in homology, as desired.
\end{proof}

By the above lemma, the ePE (co)homology of a tiling is stable under barycentric subdivision after one application, by the fact that the cell isotropy groups $\widetilde{\text{Iso}}(c,i)$ are trivial in the barycentric subdivision. Invariance under barycentric refinement allows us to deduce ePE Poincar\'{e} duality, so long as our coefficient group is suitably divisible:

\begin{theorem} \label{thm: ePE PD} Let $\mf{T}$ be an eFLC polytopal tiling. Suppose that, for some $K \in \mathbb{N}$, the coefficient ring $G$ has division by the orders of isotropy groups $\# \text{Iso}(c,K)$ for every cell $c$ of $\mf{T}$. Then we have ePE Poincar\'{e} duality $H^\bullet(\mf{T}^0;G) \cong H_{d-\bullet}(\mf{T}^0;G)$. \end{theorem}

\begin{proof} The proof is essentially identical to the proof of translational PE Poincar\'{e} duality of Theorem \ref{thm: PE PD}. All that needs to be checked is that we have invariance under refinement to the barycentric subdivision for the tiling and dual tiling, that is, that we have quasi-isomorphisms $\iota \co C_\bullet(\mf{T}^0;G) \rightarrow C_\bullet(\mf{T}^0_\Delta;G)$ and $\hat{\iota} \co C_\bullet(\hat{\mf{T}}^0;G) \rightarrow C_\bullet(\mf{T}^0_\Delta;G)$.

The cell isotropy groups $\widetilde{\text{Iso}}(c,K)$ of the tiling are quotient groups of the isotropy groups of $K$-coronas $\text{Iso}(c,K)$ by the subgroups of those transformations leaving $c$ fixed, and similarly for the dual tiling. Furthermore, any rigid motion preserving the $(K+1)$-corona of a dual cell $\hat{c}$, sending $\hat{c}$ to itself, also preserves the $K$-corona of the cell $c$ in the original tiling. It follows that the cell isotropy groups (at level $K$ for $\mf{T}$ and $K+1$ for the dual tiling) have orders which divide those of the groups $\text{Iso}(c,K)$. A unital ring which has division by $n$ also has division by any divisor of $n$, and so by Lemma \ref{lem: refinement} we have the required refinement quasi-isomorphisms $\iota$ and $\hat{\iota}$. \end{proof}

\begin{exmp} Let $\mf{T}$ be the periodic cellular tiling of $\R^2$ of unit squares whose vertices lie on the integer lattice, with the standard cellular decomposition. The cells of $\mf{T}$ have non-trivial isotropy: $\widetilde{\text{Iso}}(f,i) \cong \Z/4\Z$ for a face $f$ and $\widetilde{\text{Iso}}(e,i) \cong \Z/2\Z$ for an edge $e$. So the ePE (co)homology groups are not necessarily invariant under barycentric subdivision unless taken over coefficients $G$ with division by $4$. Since there is only one $0$-cell and one $2$-cell up to rigid motion, the ePE complexes over $\mathbb{Q}$ coefficients read
\[
0 \rightarrow \mathbb{Q} \rightarrow 0 \rightarrow \mathbb{Q} \rightarrow 0.
\]
There is no generator in degree one, since an indicator (co)chain at an edge $e$ is not invariant under the rigid motion at $e$ reversing its orientation. So the ePE invariants are $H^k(\mf{T}^0;\mathbb{Q}) \cong H_{2-k}(\mf{T}^0;\mathbb{Q}) \cong \mathbb{Q}$ for $k=0,2$ and are trivial otherwise.

To calculate over $\Z$ coefficients, we pass to the barycentric subdivision $\mf{T}_\Delta$ so that the cells have trivial isotropy. In this case we have that $H^k(\mf{T}^0_\Delta) \cong \Z$ for $k=0,2$ and are trivial otherwise. This agrees with the observation that $\Omega_\mf{T}^0$ is homeomorphic to the $2$-sphere, which by Theorem \ref{thm: Cech = ePE} has isomorphic cohomology.  For the ePE homology we have that $H_k(\mf{T}^0_\Delta) \cong \Z \oplus (\Z/2\Z) \oplus (\Z/4\Z)$, $0$, $\Z$ for $k=0,1,2$, respectively.

For this example ePE Poincar\'{e} duality $H^k(\mf{T}_\Delta^0) \cong H_{d-k}(\mf{T}_\Delta^0)$ fails over $\Z$ coefficients. Theorem \ref{thm: ePE PD} does not apply since, whilst the \emph{cells} of $\mf{T}_\Delta$ have trivial isotropy, the isotropy groups $\text{Iso}(v,i)$ of rigid motions preserving \emph{patches} are non-trivial. The ePE homology in degree zero for a periodic tiling of equilateral triangles is $\Z \oplus (\Z/2\Z) \oplus (\Z/3\Z)$. However, its associated moduli space $\Omega_\mf{T}^0$ is still the $2$-sphere,  so we see that the ePE homology is not a topological invariant of $\Omega_\mf{T}^0$ but of a finer structure. \end{exmp}

\begin{figure}
\begin{center}
\includegraphics[width=\textwidth]{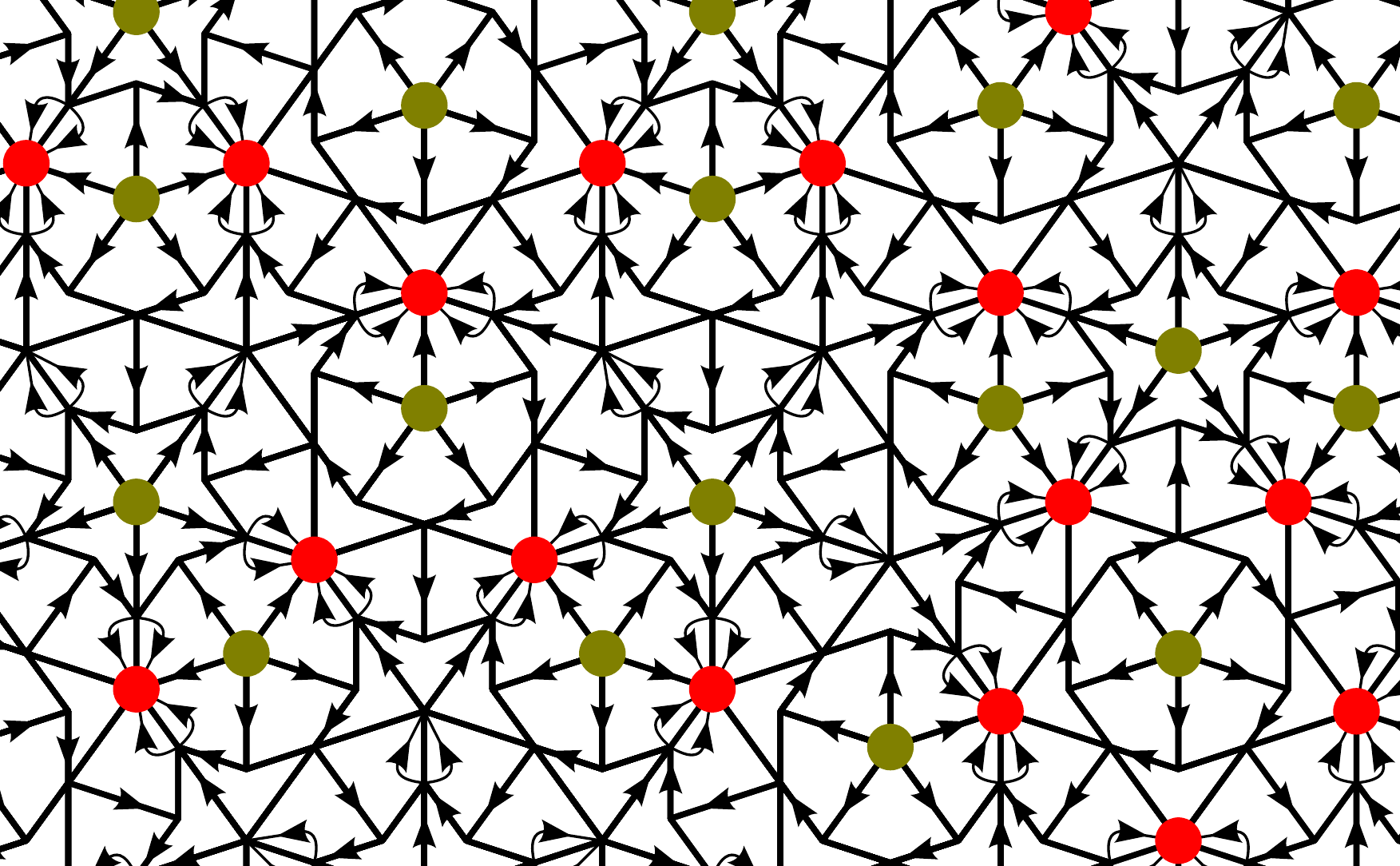} \caption{Torsion Element $t$ with $5t+\partial_1(-\cf(\text{E}1) + \cf(\text{E}2) - \cf(\text{E}4) -2 \cdot \cf(\text{E}7))=0$.}
\label{fig: S0TorBdy}
\end{center}
\end{figure}

\begin{exmp} \label{ex: Penrose ePE} Let $\mf{T}$ be a Penrose kite and dart tiling. Its ePE cohomology is
\[
H^k(\mf{T}^0) \cong \check{H}^k(\Omega_\mf{T}^0) \cong \begin{cases}
	\Z \text{ for } k=0;\\
	\Z \text{ for } k=1;\\
	\Z^2 \text{ for } k=2.\\
	\end{cases}
\]
In degrees $k=1,2$ we have ePE Poincar\'{e} duality $H_k(\mf{T}^0) \cong H^{2-k}(\mf{T}^0)$. But in degree zero, as we shall calculate in Section \ref{sec: PE Homology of Hierarchical Tilings}, we have that $H_0(\mf{T}^0) \cong \Z^2 \oplus (\Z/5\Z)$. A $0$-chain $t$ representing  a $5$-torsion homology class is depicted in Figure \ref{fig: S0TorBdy}, along with an ePE $1$-chain whose boundary is $-5t$. Specifically, the torsion element $t =\cf(\text{sun}) + \cf(\text{star}) - \cf(\text{queen})$ is a linear combination of indicator $0$-chains of certain rigid equivalence classes of star patches of $0$-cells (the seven equivalence classes of such star patches are given in Figure \ref{fig: VaETypes}). This torsion element turns out to be relevant to calculation of the \v{C}ech cohomology $\check{H}^\bullet(\Omega_\mf{T}^{\text{rot}})$ of the rigid hull \cite{Wal16}.

Generators for the free part of $H_0(\mf{T}^0)$ may be taken as $\cf(\text{sun})$ and $\cf(\text{star})$. A generator of $H_1(\mf{T}^0) \cong \Z$ is illustrated in red in Figure \ref{fig: Penrose}. It is the cycle running along the bottoms of the dart tiles, named $\rho'$ in the discussion of Example \ref{ex: Penrose PE}. \end{exmp}

\subsection{Restoring Poincar\'{e} duality} \label{subsec: Restoring PD}

The above examples show that ePE Poincar\'{e} duality can fail in the presence of non-trivial rotational symmetry. We consider the discrepancy between the ePE homology and ePE cohomology to be a feature of interest, and it will be of relevance in forthcoming work on the cohomology of Euclidean hulls \cite{Wal16}. However, to relate the ePE homology back to more familiar invariants, we shall describe here how one may modify the definition of ePE homology so as to restore duality with the ePE cohomology.

We restrict to the case that $\mf{T}$ is a tiling of $\R^2$. The higher dimensional situation is much more complicated, see the comments at the end of this subsection. We assume that $\mf{T}$ has eFLC and that it has been suitably subdivided so that any points of local rotational symmetry are contained in the vertex set of $\mf{T}$; this may be achieved for any eFLC polytopal tiling by a single barycentric subdivision.

\begin{definition} \label{def: modified complex} Define the sub-chain complex
\[
C^\dagger_\bullet(\mf{T}^0) \coloneqq 0 \leftarrow C^\dagger_0(\mf{T}^0) \xleftarrow{\partial_1} C^\dagger_1(\mf{T}^0) \xleftarrow{\partial_2} C^\dagger_2(\mf{T}^0) \leftarrow 0
\]
of the ePE complex $C_\bullet(\mf{T}^0)$ of $\mf{T}$ as follows. We let $C_k^\dagger(\mf{T}^0) \coloneqq C_k(\mf{T}^0)$ for $k = 1$, $2$. In degree zero we let $C_0^\dagger(\mf{T}^0)$ consist of those ePE chains $\sigma$ for which there exists some $i \in \N$ so that, whenever the $i$-corona of a vertex $v$ has rotational symmetry of order $n$, then the coefficient of $v$ in $\sigma$ is divisible by $n$. Denote the homology of this chain complex by $H_\bullet^\dagger(\mf{T}^0)$. \end{definition}

To see that $C_\bullet^\dagger(\mf{T}^0)$ is well-defined, firstly note that the boundary of an ePE chain is ePE, so it suffices to check that, given an ePE $1$-chain $\sigma$, there exists some $i$ for which $\partial(\sigma)$ assigns values multiples of $n$ to vertices whose $i$-coronas have $n$-fold symmetry. Since $\sigma$ is ePE, there exists some $j$ for which $\sigma$ assigns the same (oriented) coefficients to any two edges whose $j$-coronas are equivalent up to a rigid motion. Suppose that the $(j+1)$-corona of a vertex $v$ has $n$-fold rotational symmetry; these symmetries induce rigid equivalences between $j$-coronas of the edges incident with $v$. Since no edge is fixed by any non-trivial rotation, the rotations partition these edges into orbits of $n$ elements, each having equivalent $j$-coronas. It follows that the coefficient of $v$ in $\partial(\sigma)$ is some multiple of $n$, as desired.

With only minor modifications to the proof of Theorem \ref{thm: ePE PD} we obtain the following:

\begin{theorem} \label{thm: modified ePE PD} Let $\mf{T}$ be a polytopal tiling of $\R^2$ with eFLC and with points of local rotational symmetry contained in the vertex set of $\mf{T}$. Then we have Poincar\'{e} duality
\[
H^\bullet(\mf{T}^0) \cong H^\dagger_{2 - \bullet}(\mf{T}^0).
\]
\end{theorem}

\begin{proof} From the classical Poincar\'{e} duality pairing, we have an isomorphism of complexes $C^\bullet(\mf{T}^0) \cong C_{d-\bullet}(\hat{\mf{T}}^0)$ between the ePE cohomology and the ePE homology of the dual tiling. The issue with ePE Poincar\'{e} duality is that we do not necessarily have an isomorphism $H_\bullet(\hat{\mf{T}}^0) \cong H_\bullet(\mf{T}^0)$. In particular, we may not have a refinement quasi-isomorphism $\hat{\iota} \co C_\bullet(\hat{\mf{T}}^0) \rightarrow C_\bullet(\mf{T}^0_\Delta)$; the conditions of Lemma \ref{lem: refinement} are not satisfied since vertices with local rotational symmetry in $\mf{T}$ lead to dual tiles of $\hat{\mf{T}}$ with non-trivial cell isotropy.

Following the proof of Lemma \ref{lem: refinement}, we see that $\hat{\iota}$ can be made a quasi-isomorphism by replacing its range with $C^\dagger_\bullet(\mf{T}^0_\Delta)$. Let
\[
\hat{\iota}^\dagger \co C_\bullet(\hat{\mf{T}}^0) \rightarrow C^\dagger_\bullet(\mf{T}^0_\Delta)
\]
denote the canonical inclusion of chain complexes. The map $\hat{\iota}$ may fail to be a quasi-isomorphism in degree zero. It may not be the case that an ePE $0$-chain $\sigma$ of $\mf{T}_\Delta$ is homologous to a chain supported on the $0$-skeleton of $\hat{\mf{T}}$, since we are forced to `remove' $0$-chains of $C_0(\mf{T}_\Delta^0)$ from barycentres of rotationally invariant dual cells in multiples of the local symmetry at the corresponding vertices of $\mf{T}$. This issue is alleviated by passing to $C^\dagger_\bullet(\mf{T}_\Delta^0)$, since now the barycentres of such dual cells may only be assigned coefficients which are multiples of the orders of symmetries of the corresponding $i$-corona in the dual tiling for some sufficiently large $i$. The rest of the proof follows similarly to the proof of Theorem \ref{thm: ePE PD}; we end up with the following diagram of quasi-isomorphisms:
\[
C^\bullet(\mf{T}^0) \xrightarrow{- \cap \Gamma} C_{d-\bullet}(\hat{\mf{T}^0}) \xrightarrow{\hat{\iota}^\dagger} C^\dagger_{d-\bullet}(\mf{T}^0_\Delta) \xleftarrow{\iota^\dagger} C^\dagger_{d-\bullet}(\mf{T}^0).
\]\end{proof}

The above theorem shows that we may express the ePE cohomology of a two-dimensional tiling, and hence the \v{C}ech cohomology of the associated space $\Omega_\mf{T}^0$, in terms of the ePE homology of $\mf{T}$ but with certain restricted coefficients in degree zero. One may ask on the relationship between the ePE homology before and after the modification to restore Poincar\'{e} duality in degree zero. The only non-trivial chain group of the relative chain complex of the pair $C_\bullet(\mf{T}^0) \leq C_\bullet^\dagger(\mf{T}^0)$ is a torsion group in degree zero. It is not difficult to show that it is isomorphic to $\prod_{\mf{T}_i}\Z/n_i \Z$, where the product is taken over all rotation classes of tilings $\mf{T}_i$ of $\Omega_\mf{T}^0$ with $n_i$-fold rotational symmetry at the origin, at least in the case that there are only finitely many such tilings (and a similar statement holds still with infinitely many such tilings). It follows that the ePE homology and the \v{C}ech cohomology of $\Omega_\mf{T}^0$ are isomorphic in cohomological degrees one and two, and in degree zero we have that $H_0(\mf{T}^0)$ is an extension of $\check{H}^2(\Omega^0_\mf{T})$ by a torsion group determined by the rotational symmetries of tilings in the hull of $\mf{T}$.

In higher dimensions the situation is far more complicated, and we delay exposition of it to future work. The precise relationship between the ePE homology and ePE cohomology is then best expressed via a more complicated gadget, a spectral sequence analogous to the one of Zeeman \cite{Zee63}.

\begin{exmp} Consider again the periodic tiling $\mf{T}$ of unit squares. Its barycentric subdivision $\mf{T}_\Delta$ has trivial cell isotropy, but has rotational symmetry at the vertices of $\mf{T}_\Delta$ (i.e., at the barycentres of the cells of $\mf{T}$). In particular, the vertices have rotational symmetry of orders $4$, $2$ and $4$ at the vertices of $\mf{T}_\Delta$ corresponding, respectively, to the vertices, edges and faces of $\mf{T}$. So we replace the degree zero ePE chain group
\[
C_0(\mf{T}_\Delta^0) \cong \Z \langle v \rangle \oplus \Z \langle e \rangle \oplus \Z \langle f \rangle.
\]
by its modified version
\[
C_0^\dagger(\mf{T}_\Delta^0) \cong 4\Z \langle v \rangle \oplus 2\Z \langle e \rangle \oplus 4\Z \langle f \rangle.
\]
One easily computes the resulting homology group in degree zero to be $H_0^\dagger(\mf{T}_\Delta^0) \cong \Z$, restoring Poincar\'{e} duality:
\[
H_0^\dagger(\mf{T}_\Delta^0) \cong H^2(\mf{T}_\Delta^0) \cong \check{H}^2(\Omega_{\mf{T}_\Delta}^0) \cong H^2(S^2) \cong \Z.
\]
\end{exmp}

\begin{exmp} \label{ex: Penrose mePE} We saw in Example \ref{ex: Penrose ePE} that ePE Poincar\'{e} duality $H_k(\mf{T}^0) \cong H^{2-k}(\mf{T}^0)$ fails for the Penrose kite and dart tilings in homological degree $k=0$; we have extra $5$-torsion in the ePE homology, a generator is depicted in Figure \ref{fig: S0TorBdy}. Our method of calculation for the ePE homology of substitution tilings in Section \ref{sec: PE Homology of Hierarchical Tilings} may be modified to compute instead $H_0^\dagger(\mf{T}^0)$. We calculate that, indeed, Poincar\'{e} duality is restored:
\[
H_0^\dagger(\mf{T}^0) \cong H^2(\mf{T}^0) \cong \check{H}^2(\Omega_\mf{T}^0) \cong \Z^2.
\]
The modified degree zero ePE homology group is freely generated by, for example, the indicator $0$-chains of the queen and king vertex types (see Figure \ref{fig: VaETypes}). \end{exmp}

\subsection{Rotation Actions on Translational PE Cohomology} \label{subsec: Rotation Actions on Translational PE Cohomology} In the case that our tiling $\mf{T}$ is FLC (and not just eFLC) there is an alternative way of integrating the action of rotations with the PE invariants of $\mf{T}$. Firstly, we assume that some rotation group acts nicely on $\mf{T}$:

\begin{definition} We say that a finite subgroup $\Theta \leq \text{SO}(d)$ \emph{acts on $\mf{T}$ by rotations} if, for every patch $P$ of $\mf{T}$ and $g \in \Theta$, we have that $g P$ is also a patch of $\mf{T}$, up to translation. If, additionally, we have that patches $P$ and $Q$ agree up to a rigid motion only when $P$ and $g Q$ agree up to translation for some $g \in \Theta$, then we say that $\mf{T}$ \emph{has rotation group} $\Theta$. \end{definition}

A PE $k$-cochain $\psi \in C^k(\mf{T}^1)$ may be identified with a sum of indicator cochains of $i$-coronas of $k$-cells of $\mf{T}$ for some sufficiently large $i$. So if $\Theta$ acts on $\mf{T}$ by rotations, then $\Theta$ naturally acts on $C^\bullet(\mf{T}^1)$. On an indicator cochain $\cf(P_c)$, $\Theta$ acts by $g \cdot \cf(P_c) \coloneqq \cf(g^{-1} \cdot P_c)$. Replacing cellular cochains with cellular Borel--Moore chains, the same is true of the PE \emph{homology}. The action of $\Theta$ commutes with the (co)boundary maps, and so we have an action of $\Theta$ on the PE (co)homology. Furthermore, $\Theta$ naturally acts as a group of homeomorphisms on the tiling space $\Omega^1_\mf{T}$; the points of $\Omega^1_\mf{T}$ may be identified with tilings, and $\Theta$ acts by rotation at the origin. So $\Theta$ acts on the \v{C}ech cohomology $\check{H}^\bullet(\Omega^1_\mf{T})$. These actions are compatible:

\begin{proposition} \label{prop: rotation actions} Suppose that $\mf{T}$ is FLC and that $\Theta$ acts on $\mf{T}$ by rotations. Then the isomorphisms $\check{H}^\bullet(\Omega^1_\mf{T}) \cong H^\bullet(\mf{T}^1) \cong H_{d-\bullet}(\mf{T}^1)$ of theorems \ref{thm: Cech = PE} and \ref{thm: PE PD} commute with the group actions of $\Theta$. \end{proposition}

\begin{proof} The action of rotation on $\Omega^1_\mf{T}$ is canonically induced at the level of the approximants $\Gamma_i^1$. The isomorphism between the \v{C}ech cohomology of an inverse limit space and the direct limit of the cohomologies of its approximants is natural with respect to maps like this, so the isomorphisms
\[
\check{H}^\bullet(\Omega^1_\mf{T}) \cong \check{H}^\bullet(\varprojlim(\Gamma^1_i,\pi_{i,j})) \cong \varinjlim (H^\bullet(\Gamma^1_i), \pi_{i,j}^*) \cong H^\bullet(\mf{T}^1).
\]
each commute with the action of rotation. The Poincar\'{e} duality isomorphism $H^\bullet(\mf{T}^1) \cong H_{d-\bullet}(\mf{T}^1)$ of Theorem \ref{thm: PE PD} was induced by the classical pairing of a cochain with its dual chain, along with the induced maps (and their inverses) associated to barycentric refinement. Each of these maps, at the chain level, are easily seen to commute with the action of rotation.\end{proof}

Suppose that $\mf{T}$ has rotation group $\Theta$. One may ask to what extent the \v{C}ech cohomology $\check{H}^\bullet(\Omega_\mf{T}^0)$ naturally corresponds to the subgroup of $\check{H}^\bullet(\Omega^1_\mf{T})$ of elements of which are invariant under the action of $\Theta$. More concretely, we have the quotient map
\[
q \co \Omega^1_\mf{T} \rightarrow \Omega_\mf{T}^0 = \Omega^1_\mf{T} / \Theta
\]
given by identifying tilings which agree up to a rotation at the origin. Since $q = q \circ g$ for all $g \in \Theta$, the induced map
\[
q^* \co \check{H}^\bullet(\Omega_\mf{T}^0) \rightarrow \check{H}^\bullet(\Omega^1_\mf{T})
\]
has image contained in the \emph{rotationally invariant part of $\check{H}^\bullet(\Omega^1_\mf{T})$}, denoted
\[
\check{H}_\Theta^\bullet(\Omega^1_\mf{T}) \coloneqq \{[\psi] \in \check{H}^\bullet(\Omega^1_\mf{T}) \mid [\psi] = g \cdot [\psi] \text{ for all } g \in \Theta\}.
\]
Let $q^*_\Theta$ be the corestriction of $q^*$ to the range $\check{H}_\Theta^\bullet(\Omega^1_\mf{T})$. It is not difficult to show (c.f., \cite[Theorem 7]{BDHS10}) that $q^*_\Theta$ is an isomorphism when taking cohomology over divisible coefficients:

\begin{proposition} Let $\mf{T}$ be FLC with rotation group $\Theta$. If $G$ is a unital ring with division by $\# \Theta$ then $q_\Theta^* \co \check{H}^\bullet(\Omega_\mf{T}^0;G) \rightarrow \check{H}_\Theta^\bullet(\Omega^1_\mf{T};G)$ is an isomorphism.
\end{proposition}

\begin{proof} We suppose that the cell isotropy groups of $\mf{T}$ are trivial (without loss of generality, since otherwise we may simply pass to the barycentric subdivision). By the natural identification of the \v{C}ech cohomology with the PE cohomology, the map $q^*$ corresponds to the induced map of the inclusion of cochain complexes
\[
\iota \co C^\bullet(\mf{T}^0) \hookrightarrow C^\bullet(\mf{T}^1);
\]
note that $C^\bullet(\mf{T}^0)$ is the sub-cochain complex of $C^\bullet(\mf{T}^1)$ of cochains which are invariant under the action of $\Theta$. There is a self-cochain map
\[
r \co C^\bullet(\mf{T}^1) \rightarrow C^\bullet(\mf{T}^1)
\]
defined by $r(\psi) \coloneqq \sum_{g \in \Theta} g \cdot \psi$. Clearly $g \cdot r(\psi) = r(\psi)$ for any $\psi \in C^\bullet(\mf{T}^1)$, so $r$ in fact defines a map into the ePE cochain complex. We have that $r \circ \iota = \iota \circ r$ is the times $\# \Theta$ map upon restriction to the rotationally invariant part of the PE cohomology. By our assumption on the divisibility of the coefficient group $G$, this map is an isomorphism when taking cohomology over $G$ coefficients, and so $q_\Theta^*$, is an isomorphism. \end{proof}

When working over non-divisible coefficients, $q_\Theta^*$ is typically not an isomorphism. For a two-dimensional tiling, we may factor $q^*_\Theta$ through the ePE homology of $\mf{T}$:

\begin{theorem} \label{thm: rotationally invariant triangle} Suppose that $\mf{T}$ is two-dimensional, FLC, has points of local rotational symmetry contained in the vertex set of $\mf{T}$ and has rotation group $\Theta$. Then we have the following commutative triangle, with $i$ injective:
\begin{center}
\begin{tikzpicture}[auto]
\node (a1) {$\check{H}^\bullet(\Omega_\mf{T}^0)$};
\node  (a2) [below right= of a1] {$H_{2-\bullet}(\mf{T}^0)$};
\node  (a3) [above right= of a2] {$\check{H}^\bullet_\Theta(\Omega_\mf{T}^1)$};

\begin{scope}[every node/.style={scale=.7}]
\draw [right hook->] (a1) to node [swap] {$i$} (a2);
\draw [->] (a1) to node {$q_\Theta^*$} (a3);
\draw [->] (a2) to node[swap] {$f$} (a3);
\end{scope}
\end{tikzpicture}
\end{center}
\end{theorem}

\begin{proof} The inclusions of chain complexes $C_\bullet^\dagger(\mf{T}^0) \leq C_\bullet(\mf{T}^0) \leq C_\bullet(\mf{T}^1)$ induce a triangle much like the one above. The first inclusion induces an inclusion on homology, since the corresponding relative homology group is concentrated in degree zero. The rest of the proof follows from establishing that, up to isomorphism, the map $q^*$ corresponds to the induced map of the inclusion $C_\bullet^\dagger(\mf{T}^0) \leq C_\bullet(\mf{T}^1)$. This is a straightforward check, applying the ideas of the proof of Proposition \ref{prop: rotation actions}.\end{proof}

One might hope for the homomorphism $f$ in the theorem above to always be surjective. In that case, we may interpret the theorem as follows: the ePE homology extends the ePE cohomology by adding `missing' elements corresponding to PE cochains whose cohomology classes are $\Theta$-invariant, but are not actually represented by $\Theta$-invariant (ePE) cochains. Whilst $f$ usually does have a larger image than $q_\Theta^*$, the example of the Penrose tiling below shows that, in fact, $f$ need not be surjective in general.

\begin{exmp} Let $\mf{T}_\Delta$ be the barycentric subdivision of the periodic square tiling. In degree two we have that $\check{H}^2(\Omega^1_\mf{T}) \cong \Z$ which, in terms of PE cohomology, is freely generated by an indicator cochain for the square tiles, by a choice of some $2$-simplex of the barycentric subdivision of the unit square. Its cohomology class (but not the cochain itself) is invariant under rotation, so $H_\Theta^2(\Omega^1_\mf{T}) = H^2(\Omega^1_\mf{T})$.

The ePE cohomology is freely generated by the $2$-cochain which indicates each $2$-simplex of a chosen handedness, the map $q_\Theta^* \co \Z \rightarrow \Z$ is the times $4$ map in degree two. So there are classes of the PE cohomology which are invariant under rotation, but are not represented by rotationally invariant cochains. These `missing' elements are represented in the ePE homology, though. We have that $H_0(\mf{T}_\Delta) \cong \Z \oplus (\Z/2\Z) \oplus (\Z/4\Z)$ has free part generated by the $0$-chain indicating the centres of squares. In degree zero $q^*_\Theta$ factorises as $q^*_\Theta = f \circ i$ where $i$ is the $\times 4$ map onto the free component of $H_0^\dagger(\mf{T}^0)$, and the map $f$ is given by the projection $f(x,[y]_2,[z]_4) = x$.\end{exmp}

\begin{exmp} \label{ex: Penrose results} Let $\mf{T}$ be a Penrose kite and dart tiling. The cohomology $\check{H}^2(\Omega^1_\mf{T})$, along with the action of rotation by  $\Z/10\Z$ on it, is computed in \cite{Sad08}. By Proposition \ref{prop: rotation actions}, we may essentially mimic such calculations using instead PE homology. We compute (according to the method to be outlined in the next section), consistently with previous calculations, that over rational coefficients the action of rotation on $\check{H}^2(\Omega^1_\mf{T};\mathbb{Q}) \cong H_0(\mf{T};\mathbb{Q}) \cong \mathbb{Q}^8$ splits into the following irreducible sub-representations. We have two one-dimensional irreducibles corresponding to the trivial representation. There are two one-dimensional irreducibles corresponding to the representation sending the generator $[1]_{10} \in \Z/10\Z$ to the map $x \mapsto -x$. And we have a four-dimensional irreducible, the `vector representation' $\mathbb{Q}[r]/(r^4-r^3+r^2-r+1)$ which sends $[1]_{10}$ to the map $(w,x,y,z) \mapsto (-z,w+z,x-z,y+z)$.

However, over integral coefficients the representation does not decompose into irreducibles. We find elements $\sigma_1,\ldots,\sigma_4$ of $H_0(\mf{T}^1) \cong \Z^8$ upon which rotation acts trivially on $\sigma_1$ and $\sigma_2$, and sends $\sigma_3 \mapsto -\sigma_3$ and $\sigma_4 \mapsto -\sigma_4$. One may extend either of the pairs $(\sigma_1,\sigma_2)$ and $(\sigma_3,\sigma_4)$ to integer bases for $\Z^8$, but the integer span of $(\sigma_1,\sigma_2,\sigma_3,\sigma_4)$ is only an index $4$ sublattice of $H_0(\mf{T}^1)$.

In degree one, over integral coefficients, we have that the action of rotation on $\check{H}^1(\Omega^1_\mf{T}) \cong H_1(\mf{T}^1) \cong \Z^5$ decomposes to irreducibles. In terms of PE chains, it splits as the direct sum of the one-dimensional trivial representation generated by the cycle $\rho$ (see Example \ref{ex: Penrose PE}) and the four-dimensional vector representation generated by $\nu_0$ (see Example \ref{ex: Penrose PE} and Figure \ref{fig: Ammann Chain}) and its first three rotates.

We may now explain how the ePE (co)homology and the invariant part of the PE cohomology are tied together. Recall from examples \ref{ex: Penrose ePE} and \ref{ex: Penrose mePE} that
\[
\begin{array}{l}
\check{H}^1(\Omega_\mf{T}^0) \cong H^\dagger_1(\mf{T}^0) \cong H_1(\mf{T}^0) \cong \Z \langle \rho' \rangle; \\
\check{H}^2(\Omega_\mf{T}^0) \cong H^\dagger_0(\mf{T}^0) \cong \Z \langle \cf(\text{queen}) \rangle \oplus \Z \langle \cf(\text{king}) \rangle; \\
H_0(\mf{T}^0) \cong \Z \langle \cf(\text{sun}) \rangle \oplus \Z \langle \cf(\text{star}) \rangle \oplus (\Z / 5\Z) \langle \cf(\text{sun}) + \cf(\text{star}) - \cf(\text{queen}) \rangle.
\end{array}
\]
The rotationally invariant parts of the cohomology are:
\[
\check{H}^1_\Theta(\Omega^1_\mf{T}) \cong \Z; \ \ \check{H}_\Theta^2(\Omega^1_\mf{T}) \cong \Z^2.
\]
Repeating this calculation using PE homology, we find that, indeed, the rotationally invariant part of $H_0(\mf{T}^1)$ is isomorphic to $\Z^2$, and of $H_1(\mf{T}^1)$ is isomorphic to $\Z$. Furthermore, we may calculate explicit generators for these subgroups in PE homology. We find that the rotationally invariant part in degree zero is generated by PE $0$-chains $\sigma_1$, $\sigma_2$ for which $5 \sigma_1 \simeq \cf(\text{queen})$ and $5 \sigma_2 \simeq \cf(\text{king})$. The rotationally invariant part of the PE homology in degree one is freely generated by the chain $\rho$, discussed in Example \ref{ex: Penrose PE}, given by restricting the $1$-chain of Figure \ref{fig: Penrose} to loops of a chosen rotational parity.

With respect to the basis elements discussed above, we may summarise with the following commutative diagrams in cohomological degrees one and two:
\begin{center}
\begin{tikzpicture}[auto, node distance = 1cm and -0.4cm]
\node (a1)                      {$\check{H}^1(\Omega_\mf{T}^0) \cong \Z$};
\node (a2) [below right= of a1] {$H_1(\mf{T}^0) \cong \Z$};
\node (a3) [above right= of a2] {$\check{H}^1_\Theta(\Omega_\mf{T}^1) \cong \Z$};

\node (b1) [right= of a3]       {\phantom{space}};

\node (a4) [right= of b1]       {$\check{H}^2(\Omega_\mf{T}^0) \cong \Z^2$};
\node (a5) [below right= of a4] {$H_0(\mf{T}^0) \cong \Z^2 \oplus (\Z/5\Z)$};
\node (a6) [above right= of a5] {$\check{H}^2_\Theta(\Omega_\mf{T}^1) \cong \Z^2$};

\begin{scope}[every node/.style={scale=.7}]
\draw [->]            (a1) to node [swap] {$i^1 = \operatorname{id}$}   (a2);
\draw [->]            (a1) to node        {$q_\Theta^1 = \times 2$}     (a3);
\draw [->]            (a2) to node [swap] {$f^1 = \times 2$}            (a3);
\draw [->] (a4) to node[left,pos=0.6] {$
\begin{array}{c}
	i^2(x,y) = \\
	(x + 2y, x -3y,[4x+3y]_5)
\end{array}$}                                                           (a5);
\draw [->]            (a4) to node        {$q_\Theta^2(x,y) = (5x,5y)$} (a6);
\draw [->, pos = 0.7] (a5) to node[swap] {
$\begin{array}{c}
	f^2(x,y,[z]_5) = \\
	(3x+2y,x-y)
\end{array}$}                                                           (a6);
\end{scope}
\end{tikzpicture}
\end{center}

In degree one, the ePE (co)homology only corresponds to an index $2$ subgroup of the rotationally invariant part of the PE cohomology. In top cohomological degree, we have that the ePE cohomology $H^2(\mf{T}^0)$ (which corresponds to $\check{H}^2(\Omega_\mf{T}^0)$ or $H_0^\dagger(\mf{T}^0)$) maps to an index $25$ subgroup of the rotationally invariant part $\check{H}_\Theta^2(\Omega^1_\mf{T})$ of the PE cohomology. More, but not all, is added by considering instead the ePE \emph{homology}: $H_0^\dagger(\mf{T}^0) \cong H^2(\mf{T}^0)$ is an index $5$ subgroup of the ePE homology $H_0(\mf{T}^0)$, and the image of $H_0(\mf{T}^0)$ under $f$ is an index $5$ subgroup of $\check{H}_\Theta^2(\Omega^1_\mf{T})$. \end{exmp}

\subsection{Generalising The PE Framework} \label{subsec: Generalising}

Many of the constructions and results of this section did not rely on having a polytopal tiling of Euclidean space, so much as simply having a cell complex (the underlying complex of the tiling) along with a notion of when cells of that complex are equivalent to a certain radius (that is, when those cells have identical $i$-coronas, up to an agreed type of transformation). There are interesting examples of combinatorial tilings, such as the pentagonal tilings of Bowers and Stephenson \cite{BowSte97}, which are most naturally viewed as tilings of spaces which are non-Euclidean. We outline below a unified setting which allows one to deal with tilings such as this, as well as more general structures, see Example \ref{ex: solenoid}.

Recall that a CW complex is called \emph{regular} if the attaching maps of its cells may be taken to be homeomorphisms. Regular CW complexes are a sensible starting point for us here, since they allow for the construction of barycentric subdivisions and dual complexes (which, in analogy to simplicial complexes, is owing to them being essentially determined combinatorially by their face posets, see \cite{Bjo84}). Let $\mc{T}$ be a regular CW complex, it will play the r\^{o}le of the underlying cell complex of our tiling of interest.

To define the analogue of pattern-equivariance of a (co)chain, we need a notion of two (oriented) cells being equivalent in the tiling to a certain magnitude. This `magnitude' could be parametrised by, say, $ \R_{>0}$ if we want to express agreement between local patches to a certain radius, or perhaps by $\mathbb{N}$ for a combinatorial notion of patch size, such as agreement between $i$-coronas. Ultimately, there is no gain in preferring one, or indeed either of these choices. Recall that a partially ordered set $(\Lambda,\leq)$ is called \emph{directed} if for any two elements $\lambda_1, \lambda_2 \in \Lambda$ there is a third satisfying $\lambda \geq \lambda_1,\lambda_2$. We shall let some such directed set parametrise magnitude of agreement between cells of our tiling.

It is not quite enough to know which cells of our tiling are equivalent to a certain magnitude, one also needs to know \emph{how} they are equivalent. For a cell $c$ of $\mc{T}$, its (\emph{closed}) \emph{star} $\text{St}_\mc{T}(c)$ is defined to be the subcomplex of $\mc{T}$ whose support is the set of cells containing $c$. The way in which a cell $c_1$ may be considered as `equivalent' to a cell $c_2$ to some magnitude $\lambda$ will be recorded by a finite set of cellular homeomorphisms $\Phi \co \text{St}_\mc{T}(c_1) \rightarrow \text{St}_\mc{T}(c_2)$. The star of a cell defines a neighbourhood of that cell, so we may consider such morphisms as defining germs of maps by which two cells are equivalent.

With these interpretations of the ingredients, we may give the following definition, which provides a structure upon which one may define PE (co)homology and various other related constructions. The axioms will be further motivated below.

\begin{definition}\label{def: SIS} A \emph{system of internal symmetries} (or \emph{SIS}, for short) $\mf{T}$ consists of the following data:
\begin{itemize}
	\item A finite-dimensional and locally-finite regular CW complex $\mc{T}$.
	\item A directed set $(\Lambda,\leq)$ called the \emph{magnitude poset}.
	\item For each $\lambda \in \Lambda$ and each pair of cells $a,b \in \mc{T}$ a set $\mf{T}_{a,b}^\lambda$ of cellular homeomorphisms $\Phi \co \text{St}_\mc{T}(a) \to \text{St}_\mc{T}(b)$ sending $a$ to $b$. We denote the collection of all such morphisms by $\mf{T}^\lambda$.
\end{itemize}

This data is required to satisfy the following:
\begin{enumerate}
	\item[(G1)] For all $\lambda \in \Lambda$ and $a \in \mc{T}$ we have that $\operatorname{id}_{\text{St}(a)} \in \mf{T}_{a,a}^\lambda$.
		\item[(G2)] For all $\lambda \in \Lambda$ and $\Phi \in \mf{T}_{a,b}^\lambda$ we have that $\Phi^{-1} \in \mf{T}_{b,a}^\lambda$.
	\item[(G3)] For all $\lambda \in \Lambda$, $\Phi_1 \in \mf{T}_{a,b}^\lambda$ and $\Phi_2 \in \mf{T}_{b,c}^\lambda$, we have that $\Phi_2 \circ \Phi_1 \in \mf{T}_{a,c}^\lambda$.
	\item[(Inc)] For all $\lambda_1 \leq \lambda_2$ we have that $\mf{T}^{\lambda_1} \supseteq \mf{T}^{\lambda_2}$.
	\item[(Res)] For all $\lambda \in \Lambda$ there exists some $\lambda_{\text{res}} \geq \lambda$ satisfying the following. Given any $b \in \mc{T}$ and face $a \preceq b$, every morphism of $\mf{T}^{\lambda_{\text{res}}}_{a,-}$ restricts to a morphism of $\mf{T}^\lambda_{b,-}$.
	\item[($\widehat{\text{Res}}$)] Dually, for all $\lambda \in \Lambda$ there exists some $\lambda_{\widehat{\text{res}}}$ satisfying the following. Given any $a \in \mc{T}$ and coface $b \succeq a$, every morphism of $\mf{T}_{b,-}^{\lambda_{\widehat{\text{res}}}}$ is a restriction of some morphism of $\mf{T}_{a,-}^\lambda$.
\end{enumerate}
\end{definition}

As explained above, morphisms $\Phi \in \mf{T}_{a,b}^\lambda$ should be interpreted as recording that \emph{cell $a$ is equivalent to cell $b$ to magnitude $\lambda$ via $\Phi$}. The groupoid axioms (G1)--(G3) state that such morphisms should include the identity morphisms, be invertible and that compatible morphisms are composable. The inclusion axiom (Inc) simply states that if two cells are 	equivalent to magnitude $\lambda_2$, via morphism $\Phi$, then they are still equivalent via $\Phi$ to any smaller magnitude $\lambda_1 \leq \lambda_2$. In short, for $\lambda_1 \leq \lambda_2$ we have an inclusion of groupoids $\mf{T}^{\lambda_2} \hookrightarrow \mf{T}^{\lambda_1}$. The final two axioms (Res) and ($\widehat{\text{Res}}$) of restriction and corestriction establish a coherence between the cellular structure of $\mc{T}$ and the restrictions between the various morphisms of $\mf{T}$.

The results of this section may be generalised to systems of internal symmetries with only minor modifications to the definitions and proofs, although these proofs are most efficiently given in a combinatorial setting. We may derive direct analogues of the following for SISs:

\begin{enumerate}
	\item PE or ePE (co)homology.
	\item The tiling space inverse limit presentation $\Omega_\mf{T} = \varprojlim_{\lambda \in \Lambda} (\Gamma_\lambda,\pi_{\lambda,\mu})$ where the approximants $\Gamma_\lambda$ are CW complexes given by identifying cells of $\mc{T}$ which are equivalent to magnitude $\lambda$.
	\item Theorem \ref{thm: Cech = ePE}, that the \v{C}ech cohomology of $\Omega_\mf{T}$ agrees with the PE cohomology of $\mf{T}$ when the cells of $\mf{T}$ (for sufficiently large magnitudes) have trivial isotropy.
	\item Lemma \ref{lem: refinement}, that we have invariance of the PE (co)homology of $\mf{T}$ over $G$ coefficients under barycentric refinement whenever $G$ (for sufficiently large magnitudes) has divisibility by the order of isotropy of cells in $\mf{T}$.
	\item Theorem \ref{thm: ePE PD}, ePE Poincar\'{e} duality $H^\bullet(\mf{T};G) \cong H_{d-\bullet}(\mf{T};G)$. For this to hold, we need firstly that the ambient space of $\mf{T}$ is a $G$-oriented $d$-manifold (or even just homology $G$-manifold) with pattern-equivariant fundamental class $\Gamma$. We also require the analogous condition on the divisibility of the coefficient ring; that is, there exists some $\lambda \in \Lambda$ for which $G$ has division by the order of isotropy groups $\mf{T}_{a,a}^\lambda$ at every cell $a$ of $\mf{T}$.	
\end{enumerate}

\begin{exmp} \label{ex: BSP} The initial insight in Penrose's discovery of his famous tilings was that ``a regular pentagon can be subdivided into six smaller ones, leaving only five slim triangular gaps'' \cite{Pen84}. In \cite{BowSte97} Bowers and Stephenson took a similar subdivision but chose, instead of methodically filling in the slim triangular gaps, to simply remove them by identifying edges of the pentagons. Of course, this cannot be achieved in Euclidean space with regular pentagons; the result is a \emph{combinatorial} substitution. One may produce, in an analogous way to in the Euclidean setting of \cite{AndPut98,Fra08}, limiting combinatorial tilings. Declaring that each $2$-cell of such a combinatorial tiling should metrically correspond to a regular pentagon, the resulting tilings are of spaces which are homeomorphic, but not isometric to Euclidean $2$-space.

There is no notion of translation on the ambient spaces of these tilings, but there is of orientation. Let $\mc{T}$ be a Bowers--Stephenson pentagonal tiling, which we consider here simply as a regular cell complex with a choice of identification of each $2$-cell with the regular pentagon. We may define a corresponding SIS $\mf{T}^0$ as follows. Given cells $a,b \in \mc{T}$ and $i \in \N$, consider the collection of maps taking the $i$-corona of $a$ to the $i$-corona of $b$, preserving orientation and distances on each pentagonal tile (such a map is, of course, determined combinatorially by how it acts on cells). We let $(\mf{T}^0)_{a,b}^i$ be the set of such maps restricted to the stars of $a$ and $b$.

We may construct from this data an inverse limit of approximants and associated inverse limit space $\Omega_\mf{T}^0$, analogously to the Euclidean setting. The points of $\Omega_\mf{T}^0$ may be identified with pointed Bowers--Stephenson pentagonal tilings, two being `close' if they agree via an orientation-preserving isometry on large patches about their origins, up to a small perturbation of their origins. We have analogues of the ePE (co)homology groups, and of theorems \ref{thm: ePE PD} and \ref{thm: modified ePE PD}. We may not identify the (integer) ePE cohomology with the \v{C}ech cohomology of $\Omega_\mf{T}^0$, since the cells have non-trivial isotropy, but we may after replacing the tiling with its barycentric subdivision. These cohomology groups are Poincar\'{e} dual to a modification of the ePE homology of $\mf{T}_\Delta^0$, defined in an analogous fashion to the modified chain complexes $C_\bullet^\dagger(\mf{T}_\Delta^0)$ of Definition \ref{def: modified complex}.

Our method of computation of the ePE homology groups in Section \ref{sec: PE Homology of Hierarchical Tilings} easily generalises to examples such as this, we find that
\[
\begin{array}{l}
\check{H}^1(\Omega_\mf{T}^0) \cong H^\dagger_1(\mf{T}_\Delta^0) \cong H_1(\mf{T}_\Delta^0) \cong 0; \\
\check{H}^2(\Omega_\mf{T}^0) \cong H^\dagger_0(\mf{T}_\Delta^0) \cong \Z \oplus \Z[1/6]; \\
H_0(\mf{T}_\Delta^0) \cong \Z \oplus \Z[1/6].
\end{array}
\]
\end{exmp}

\begin{exmp} \label{ex: solenoid} In this example we shall see how more general objects are also naturally captured in this framework. The magnitude poset will be $\mathbb{N}$, but endowed with the partial ordering of $m \leq n$ if $m$ divides $n$; note that $(\mathbb{N},\divides)$ is a directed set. Let $\mc{T}$ be the standard cellular decomposition of $\R^d$ associated to the tiling of unit cubes with vertices at the lattice points $\Z^d$. If cells $a,b \in \mc{T}$ are equal up to a translation in $n \Z^d$, then we let $(\mf{T}^1)_{a,b}^n$ consist of the single map given by the restriction of this translation between the stars of $a$ and $b$. Otherwise, we set $(\mf{T}^1)_{a,b}^n = \emptyset$.

It is easy to check that $\mf{T}^1$ thus defined satisfies the conditions of an SIS. A (co)chain of $\mc{T}$ is PE with respect to $\mf{T}^1$ if and only if it is invariant under translation by some full rank sub-lattice of $\Z^d$. We have trivial isotropy (everything is generated by translations), and the analogous theorems and constructions of the previous sections apply. For example, in dimension $d=1$ we may calculate that
\[
\begin{array}{l}
\check{H}^0(\Omega_\mf{T}^1) \cong H^0(\mf{T}^1) \cong H_1(\mf{T}^1) \cong \Z; \\
\check{H}^1(\Omega_\mf{T}^1) \cong H^1(\mf{T}^1) \cong H_0(\mf{T}^1) \cong \mathbb{Q}.
\end{array}
\]
The first isomorphisms are given by the analogue of Theorem \ref{thm: Cech = PE} and the second by the analogue of PE Poincar\'{e} duality of Theorem \ref{thm: PE PD}. The tiling space $\Omega_\mf{T}^1$ is homeomorphic---in a fashion entirely analogous to the G\"{a}hler construction---to the inverse limit
\[
\Omega_\mf{T}^1 = \varprojlim (S^1,\pi_{m,n})
\]
where, for $m \divides n$, the map $\pi_{m,n}$ is the standard degree $(n/m)$ covering map of $S^1$. The degree one PE homology group $H_1(\mf{T}^1)$ is generated by a Borel--Moore fundamental class for $\R^1$ and the homology class $p/q \in H_0(\mf{T}^1) \cong \mathbb{Q}$ is represented by, for example, the Borel--Moore $0$-chain which assigns value $p$ to each of the vertices of $q \Z$. Note that the sequence $n_i \coloneqq i!$ is linearly ordered and cofinal in $(\mathbb{N},\divides)$, so the tiling space could instead be expressed as
\[
\Omega_\mf{T}^1 = \varprojlim (S^1 \xleftarrow{\times 2} S^1 \xleftarrow{\times 3} S^1 \xleftarrow{\times 4} S^1 \xleftarrow{\times 5} \cdots).
\]

We may restrict the construction of $\mf{T}^1$ above to the linearly ordered subset of magnitudes $\{2^n \mid n \in \mathbb{N}_0\}$. In this case, a (co)chain is PE if and only if it invariant under translation by $2^n \Z^d$ for some $n \in \mathbb{N}_0$. For $d=1$, the corresponding tiling space is the dyadic solenoid
\[
\Omega_\mf{T}^1 = \varprojlim (S^1 \xleftarrow{\times 2} S^1 \xleftarrow{\times 2} S^1 \xleftarrow{\times 2} S^1 \xleftarrow{\times 2} \cdots).
\]
Whilst these examples are not given by tilings of finite local complexity, they are close in spirit. Indeed, one may think of the dyadic example above as a hierarchical tiling. The tiles (unit cubes) may be grouped into supertiles (cubes of side-length $2$), which may be grouped into level $n$ supertiles (cubes of side-length $2^n$). However, the groupings of tiles cannot be determined using local geometric information in the underlying tiling; one says that the substitution corresponding to this example is \emph{non-recognisable}. The resulting system of internal symmetries is what one would get if the tiling were capable of deducing such an imposed hierarchy from local geometric information. For an alternative derivation of the dyadic solenoid as the tiling space of an infinite local complexity tiling, see \cite{PrFSad14ILC}.\end{exmp}


\section{PE Homology of Hierarchical Tilings} \label{sec: PE Homology of Hierarchical Tilings}

\subsection{Substitution Tilings and Their Hulls} The two main approaches to producing interesting examples of aperiodic tilings, such as the Penrose tilings, are through the cut-and-project method (see \cite{FHK02}) and through tiling substitutions (see \cite{AndPut98}). A substitution rule consists of a finite collection of \emph{prototiles} of $\R^d$, a rule for subdividing them and an expanding dilation which, when applied to the subdivided prototiles, defines patches of translates of the original prototiles. By iterating the substitution and inflating, one produces successively larger patches. A tiling is said to be \emph{admitted by the substitution rule} if every finite patch of it is a sub-patch of a translate of some iteratively substituted prototile. In analogy with symbolic dynamics, one may think of the substitution rule as generating the allowed language for a family of tilings.

Under certain conditions on the substitution rule, tilings admitted by it exist and, in addition, for each such tiling $\mf{T}_0$ there is a \emph{supertiling} $\mf{T}_1$, based on inflated versions of the prototiles, which subdivides to $\mf{T}_0$ and is itself a (dilation of an) admitted tiling. So $\mf{T}_0$ has a hierarchical structure: there is an infinite list of substitution tilings $\mf{T}_0$, $\mf{T}_1$, $\mf{T}_2$, \ldots of progressively larger tiles for which the tiles of $\mf{T}_n$ may be grouped to form the tiles of $\mf{T}_{n+1}$, with the substitution decomposing $\mf{T}_{n+1}$ to $\mf{T}_n$. For fuller details on the definition of substitution rules and their tilings, we refer the reader to \cite{AndPut98, Fra08, Sad08}.

To compute the \v{C}ech cohomology of a substitution tiling space $\Omega^1$, one typically constructs a finite CW complex $\Gamma$ along with a self-map $f$ of $\Gamma$ for which
\[
\Omega^1 \cong \varprojlim (\Gamma \xleftarrow{f} \Gamma \xleftarrow{f} \Gamma \xleftarrow{f} \cdots).
\]
The CW complex $\Gamma$ may be defined in terms of the short-range combinatorics of the patches of the substitution tilings, and the map $f$ by the action of substitution. This makes the \v{C}ech cohomology of a substitution tiling space computable.

Anderson and Putnam showed that when the substitution rule has a property known as \emph{forcing the border}, one may take $\Gamma$ as, what is now known as, the \emph{AP complex} \cite{AndPut98}, which is precisely the level zero G\"{a}hler complex $\Gamma^1_0$ (see Subsection \ref{sec: Inverse Limit Presentations}). If the substitution fails to force the border, one may work with the \emph{collared AP complex} $\Gamma^1_1$ instead. Whilst conceptually simple, passing to the collared complex can be computationally demanding; even for relatively simple substitution rules, the number of collared tiles can be unwieldy. A powerful alternative approach was developed by Barge, Diamond, Hunton and Sadun \cite{BDHS10}, which typically results in much smaller cochain complexes than for the collared AP complex. One constructs a CW complex $K_\epsilon$ by, instead of collaring tiles, collaring \emph{points} of the ambient space of the tiling. A point of $K_\epsilon$ is a description of how to tile an $\epsilon$-neighbourhood of the origin. The self map on the complex defined by the substitution is not cellular, but for small $\epsilon$ may be made homotopic to a cellular map in a canonical way, which is sufficient for cohomology computations. Another advantage of the BDHS approach is that the resulting inverse system possesses natural stratifications, which are useful in breaking down the calculations to something more tractable.

\subsection{Overview of PE Homology Approach to Calculation} We shall present below a method of calculation of the PE homology of a substitution tiling. There are various motivations for introducing it. Firstly, as we shall see, the `approximant complexes' and `connecting map' of the method are constructed from the combinatorial information of the substitution in a very direct way, which makes the approach highly amenable to computer implementation. An early implementation has been coded by the author, in collaboration with James Cranch, in the programming language Haskell, and at present is applicable to any polytopal substitution tiling of arbitrary dimension (although, aside from cubical substitutions, efficiently communicating the combinatorial information of a tiling substitution to the program is still problematic, an issue shared with any machine computation of the cohomology of substitution tilings). The approximant chain complexes of the method are much smaller than for the collared AP method; the combinatorial information required from the patches is the same as that for the BDHS approach.

Another reason for introducing this method is that it may be used to find the \emph{ePE} homology of a substitution tiling (see Example \ref{ex: Penrose ePE Calculation}) which, as we have seen, yields different information to the cohomology calculations. Furthermore, the method provides explicit generators in terms of pattern-equivariant chains. The result for the ePE homology of the Penrose tiling, along with precise descriptions of the generators of the ePE homology, will be essential in forthcoming work \cite{Wal16}.

In the translational setting, our approach can be seen to produce isomorphic direct limit diagrams to approximant homologies of the BDHS method \cite{BDHS10}, at least after collaring points of the tiling for the BDHS approximants in a way compatible with the combinatorics of the tiling (although the method that we shall describe provides a more combinatorial way of determining this diagram). The argument proceeds via a stratification of the BDHS approximants (although one which is not preserved by the connecting maps) or by applying a certain homotopy to the projective system of BDHS approximants. However, the full details of this seem to be technical, at least in general dimensions, and we avoid providing them here. In any case, the approximant complexes and connecting maps between the approximant homologies that we shall define are most naturally described in the PE homology framework. The approximant complexes used here are precisely the duals of those used by Gon\c{c}alves \cite{Gon11} in his computation of the $K$-theory of the $C^*$-algebra associated to the stable relation of a one or two-dimensional substitution tiling. This $K$-theory appears to be dual in a certain sense to the $K$-theory of the hull (that is, of the unstable relation). The fact that our technique---which involves the duals of the approximant complexes of \cite{Gon11}---calculates the (regraded) \v{C}ech cohomology of the hull of a substitution tiling seems to confirm this duality. A complete confirmation of the relationship would, however, require consideration of the connecting maps of each method of calculation.

\subsection{The Method of Computation} 

\subsubsection{The Approximant Complex} We shall assume throughout that $\mf{T} = \mf{T}_0$ is a cellular FLC tiling admitted by a primitive, recognisable, polytopal substitution rule $\omega$ with inflation constant $\lambda > 1$. We refer the reader to \cite{AndPut98} for the notion of a primitive polytopal substitution. We note that many of these assumptions may be relaxed substantially; instead of letting that detain us here, we shall discuss various generalisations in Subsection \ref{subsec: Generalisations}.

Being a \emph{recognisable} substitution means that for any tiling $\mf{T}_0$ admitted by $\omega$ there exists a unique FLC tiling $\mf{T}_1$, based on prototiles which are $\lambda$-inflations of the original prototiles, for which
\begin{itemize}
	\item the rescaled tiling $\lambda^{-1} (\mf{T}_1)$ is also admitted by $\omega$;
	\item $\omega(\mf{T}_1) = \mf{T}_0$;
	\item $\mf{T}_1$ is locally derivable from $\mf{T}_0$.
\end{itemize}
The first item simply states that the \emph{supertiling} $\mf{T}_1$ is itself an inflate of an admitted tiling. The second says that the substitution rule decomposes the supertiles of $\mf{T}_1$ to the tiles of $\mf{T}_0$. Thinking upside-down, one may group the tiles of $\mf{T}_0$ so as to form the supertiling $\mf{T}_1$. The third item states that this grouping may be performed using only local information. Since $\omega(\mf{T}_1) = \mf{T}_0$ implies that $\mf{T}_0$ is locally derivable from $\mf{T}_1$, the two are MLD. This process may be repeated, yielding a hierarchy of tilings $\{\mf{T}_n\}_{n \in \N_0}$. Each $\mf{T}_n$ carries a polytopal decomposition $\mc{T}_n$, with $\mc{T}_i$ refining $\mc{T}_j$ for $i \leq j$.

Given a $k$-cell $c$ of $\mc{T}_0$, we name the pair of $c$ along with the set of tiles properly containing $c$ the \emph{star} of $c$. Henceforth, the translation class of such a star (where translations preserve labels, if the cells are labelled) will be simply called a \emph{star}, or a \emph{$k$-star} if we wish to specify the dimension of the \emph{central cell} $c$. The first step of the calculation is to enumerate the set of stars. This may be efficiently performed algorithmically as follows:
\begin{enumerate}
	\item Begin with the set of $d$-stars, consisting of the prototiles of $\mf{T}$. Put these stars into sets $S^{\text{new}}_0$ and $S^{\text{acc}}_0$.
	\item Suppose that $S^{\text{new}}_n \neq \emptyset$ and $S^{\text{acc}}_n$ have been constructed. Substitute each star of $S^{\text{new}}_n$ and find all stars whose centres are contained in the substituted central (open) cell of the original star. All such stars which are not already elements of $S^{\text{acc}}_n$ define the set $S_{n+1}^{\text{new}}$, and are added to $S^{\text{acc}}_n$ to define $S^{\text{acc}}_{n+1}$.
	\item If $S^{\text{new}}_n$ and $S_n^{\text{acc}}$ have been constructed but $S^{\text{new}}_n = \emptyset$, then the process is terminated and the full list of stars is $S \coloneqq S^{\text{acc}}_n$.
\end{enumerate}

The stars, and the incidences between them, define our approximant chain complex:

\begin{definition} \label{def: approximant chain complex}
We define the \emph{approximant chain complex}
\[
C_\bullet^{(0)}(\mf{T}^1) \coloneqq 0 \leftarrow C_0^{(0)}(\mf{T}^1) \xleftarrow{\partial_1} C_1^{(0)}(\mf{T}^1) \xleftarrow{\partial_2} \cdots \xleftarrow{\partial_d} C_d^{(0)}(\mf{T}^1) \leftarrow 0
\]
as follows. The degree $k$ chain group $C_k^{(0)}(\mf{T}^1)$ is freely generated by the $k$-stars; $C_k^{(0)}(\mf{T}^1) \cong \Z^n$, where $n$ is the number $k$-stars. The star of a $(k-1)$-cell $c'$ determines the star of any $k$-cell $c$ containing $c'$. Orienting the central cell of each star, the boundary maps $\partial_k$ are induced from the standard cellular boundary maps by
\[
\partial_k(s) = \sum_{(k-1)\text{-stars } s'} [s',s] \cdot s'
\]
and extending linearly. Here, $[s',s]$ is the incidence number between the $(k-1)$-star $s'$ and $k$-star $s$. It is defined, for fixed $s'$ and $s$, to be the sum of incidence numbers $[c',c]$ where $c'$ is the central $(k-1)$-cell of $s'$ and $c$ is a $k$-cell of $s'$ whose star is $s$. The homology of the approximant chain complexes is the \emph{approximant homology} $H^{(0)}_\bullet(\mf{T}^1)$.
\end{definition}

Note that since we are considering translation classes of stars, it may be that there are multiple occurrences of a $k$-star $s$ in a $(k-1)$-star $s'$. An instructive perspective on the definition of the approximant chain complex is to identify a generator star $s$ with the PE indicator chain $\cf(s) \in C_k(\mf{T}^1)$, the $k$-chain given by the (infinite) sum of $k$-cells in $\mf{T} = \mf{T}_0$ which are the centres of $s$ in the tiling. With this identification, the boundary maps of the approximant complex correspond to the standard cellular boundary maps of $C_\bullet(\mf{T}^1)$ defined in the ambient tiling. That is, $C_\bullet^{(0)}(\mf{T}^1)$ is the sub-chain complex of the PE complex $C_\bullet(\mf{T}^1)$ consisting of those chains which, at any given cell $c$, depend only on the tiles of $\mf{T}$ properly containing $c$.

\begin{figure}
\begin{center}
\includegraphics[scale=0.3]{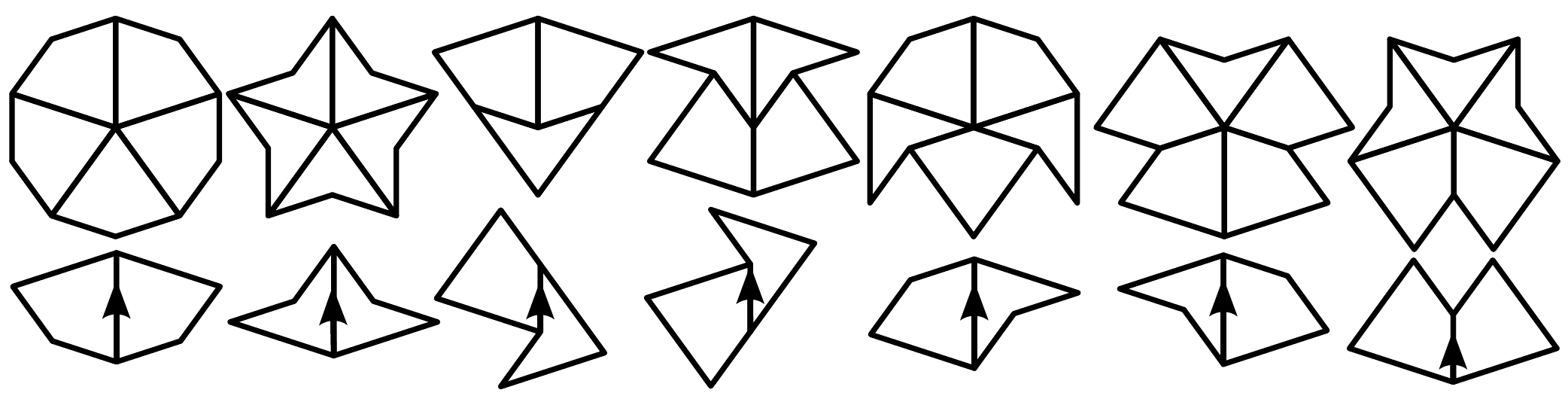}
\caption{The seven rigid equivalence classes of $0$-stars types named, in Conway's notation, `sun', `star', `ace', `deuce', `jack', `queen' and `king', in this respective order, followed by the seven $1$-stars $\text{E}1$--$\text{E}7$.}
\label{fig: VaETypes}
\end{center}
\end{figure}

\begin{exmp} \label{ex: 1d approximants} In the one-dimensional case, we may identify our tiling with a bi-infinite sequence $s \in \mc{A}^\Z$ over a finite alphabet $\mc{A}$, and the substitution rule with a map $\omega \colon \mc{A} \rightarrow \mc{A}^*$ from the alphabet $\mc{A}$ to the set of non-empty words in $\mc{A}$. A $0$-star then corresponds to a two-letter word from $\mc{A}^2$ which appears in $s$; let us denote the set of such \emph{admissible} two-letter words by $\mc{A}^2_\omega$. A $1$-star is simply a tile type, an element of $\mc{A}$. So the approximant complex is given by
\[
0 \leftarrow \Z^m \xleftarrow{\partial_1} \Z^n \leftarrow 0
\]
where $m = \# \mc{A}_\omega^2$ and $n = \# \mc{A}$. The boundary map $\partial_1$ is defined on $a \in \mc{A}$ as the formal sum of admissible two-letter words whose first letter is $a$ minus those whose second letter is $a$.\end{exmp}

\begin{exmp} One may easily verify that, up to rigid motion, there are seven distinct ways for a patch of Penrose kite and dart tiles to meet at a vertex, and seven ways for them to meet at an edge. These stars are given in Figure \ref{fig: VaETypes}. Any rotate of such a patch by some $2\pi k /10$ appears in a Penrose kite and dart tiling, so there are $54$ $0$-stars (the `sun' and `star' vertices are preserved by rotation by $2\pi /5$), there are $70$ $1$-stars and there are $20$ $2$-stars, corresponding to the rotates of the kite and dart tiles. So the approximant chain complex is of the form
\[
C_\bullet^{(0)}(\mf{T}^1) = 0 \leftarrow \Z^{54} \xleftarrow{\partial_1} \Z^{70} \xleftarrow{\partial_2} \Z^{20} \leftarrow 0.
\]
The boundary maps have a simple description in terms of the standard cellular boundary maps. For example, the boundary of the indicator chain of the $\text{E}1$ edge is given by
\[
\partial_1(\cf(r^0 \text{E}1)) = \cf(r^1 \text{sun}) +\cf(r^5 \text{jack}) - \cf(r^0 \text{ace}) - \cf(r^5 \text{queen});
\]
the head of an $\text{E}1$ edge is always a `sun' or `jack' vertex, and the tail is always an `ace' or `queen'. The notation $r^k$ above indicates that the named patch has been rotated by $2\pi k /10$ relative to its depiction in Figure \ref{fig: VaETypes}. We compute the approximant homology groups as
\[
H^{(0)}_k(\mf{T}^1) \cong  \begin{cases}
	\Z^8 \text{ for } k=0;\\
	\Z^5 \text{ for } k=1;\\
	\Z \text{ for } k=2.\\
	\end{cases}
\]
Representative cycles for the generators of these approximant homology groups are as in Example \ref{ex: Penrose PE}. \end{exmp}

\subsubsection{The Connecting Map}

In the case of the Penrose kite and dart substitution, the inclusion of the approximant chain complex into the full PE chain complex is a quasi-isomorphism. This is rarely true in general. We shall now describe how one constructs a homomorphism
\[
f \co H^{(0)}_\bullet(\mf{T}^1) \rightarrow H^{(0)}_\bullet(\mf{T}^1)
\]
from the approximant homology to itself, called the \emph{connecting map}, for which the PE homology $H_\bullet(\mf{T}^1)$ is isomorphic to the direct limit $\varinjlim (H_\bullet^{(0)}(\mf{T}^1),f)$.

To construct $f$, we now need to consider the passage from the tiling $\mf{T}_0$ to its supertiling $\mf{T}_1$. By the explanation following Definition \ref{def: approximant chain complex}, we may identify $C_\bullet^{(0)}(\mf{T}^1)$ with the sub-chain complex of PE chains which are determined at a cell $c$ by the patch of tiles properly containing $c$. To define $f$, we firstly define an auxiliary chain complex $C_\bullet^{(1)}(\mf{T}^1)$, still a subcomplex of $C_\bullet(\mf{T}^1)$, which consists of those PE chains of $\mf{T} = \mf{T}_0$ which are determined at any given cell of $\mc{T}_0$ by only which \emph{super}tiles of $\mf{T}_1$ properly contain $c$. An alternative take on this is that one has a new set of generators which, in degree $k$, are given by $k$-cells of the substituted central cells of the original stars. The boundary maps of $C_\bullet^{(1)}(\mf{T}^1)$ are induced from the standard cellular boundary maps, analogously to $C_\bullet^{(0)}(\mf{T}^1)$. That is, for a translation class of oriented $k$-cell $s$ from $\mc{T}_0$, labelled by the patch information of supertiles which contain it, we define $\partial_k(s)$ by identifying $s$ with the PE indicator chain $\cf(s) \in C_k(\mf{T}^1)$, and then define $\partial_k(s)$ to be the element of $C_{k-1}^{(1)}(\mf{T}^1)$ corresponding to $\partial(\cf(s)) \in C_{k-1}(\mf{T}^1)$.

There are two very natural maps from $C_\bullet^{(0)}(\mf{T}^1)$ to $C_\bullet^{(1)}(\mf{T}^1)$ which we shall use to construct $f$. Let $s$ be some $k$-star, which we identify with the PE indicator chain $\cf(s) \in C_k(\mf{T}^1)$. Since $\cf(s)$ is determined at any cell by the patch of supertiles properly containing that cell, $\cf(s)$ is also an element of $C_k^{(1)}(\mf{T}^1)$. We define the chain map $\iota$ as this inclusion of chain complexes.

The combinatorics, and hence stars, of $\mf{T}_0$ and $\mf{T}_1$ are identical. Since we may identify the stars of each, we may associate any chain $\sigma \in C_k^{(0)}(\mf{T}^1)$ with a chain $\sigma' \in C_k^{\text{BM}}(\mc{T}_1)$ which assigns coefficients to $k$-cells of $\mc{T}_1$ based on their neighbourhood stars in $\mf{T}_1$ identically to how $\sigma$ does in $\mf{T}_0$. The chain map $q$ is given by identifying $\sigma'$ with its representation on the finer subcomplex $\mc{T}_0$, induced by identifying the elementary chain of a $k$-cell $c$ of $\mc{T}_1$ with the sum of $k$-cells of $\mc{T}_0$ contained in $c$, suitably oriented with respect to $c$. It is easily seen that $q(\sigma) \in C_k^{(1)}(\mf{T}^1)$. The chain map $q$ is in some sense simply a refinement. So as one may expect, $q$ is a quasi-isomorphism; that is, the induced map on homology
\[
q_* \co H_\bullet^{(0)}(\mf{T}^1) \rightarrow H_\bullet^{(1)}(\mf{T}^1)
\]
is an isomorphism, see Lemma \ref{lem: hierarchical refinement}.

\begin{definition}
We define the \emph{connecting map} as
\[
f \coloneqq (q_*)^{-1} \circ \iota_* \co H_\bullet^{(0)}(\mf{T}^1) \rightarrow H_\bullet^{(0)}(\mf{T}^1),
\]
where the chain maps $\iota$ and $q$ are defined as above.
\end{definition}

\begin{theorem} \label{thm: method} There is a canonical isomorphism $\varinjlim (H_\bullet^{(0)}(\mf{T}^1),f) \cong H_\bullet(\mf{T}^1)$. \end{theorem}

By \emph{canonical} here, we mean that the isomorphism is induced by a natural association of cycles of the direct limit with PE cycles of $C_\bullet(\mf{T}^1)$; a chain at the $n$th level of the direct limit corresponds to a PE chain which only depends cell-wise on its immediate surroundings in the level $n$ supertiling $\mf{T}_n$. We delay the details of the proof to the final subsection.

We remark that the connecting map $f$ is \emph{not} canonically induced by a chain map from $C_\bullet^{(0)}(\mf{T}^1)$ to itself; it is defined at the level of homology rather than of chains. Nonetheless, we may provide a geometrically intuitive picture of the action of $f$ on the approximant cycles. Suppose that $\sigma \in C_k^{(0)}(\mf{T}^1)$ is a cycle, which we may identify with a PE $k$-cycle of $\mf{T}_0$ that only depends at any $k$-cell $c$ on the patch of tiles containing $c$. Considered as a chain living inside the supertiling $\mf{T}_1$ (but still a chain of the finer complex $\mc{T}_0$), $\sigma$ may no longer be supported on the $k$-skeleton of $\mf{T}_1$, but it is still determined cell-wise by the local patches in $\mf{T}_1$. Due to the homological properties of the cells, we may find a chain $\tau \in C_{k+1}^{(1)}(\mf{T}^1)$ for which $\sigma' \coloneqq \sigma + \partial(\tau)$ is supported on the $k$-skeleton of $\mf{T}_1$. The combinatorics of $\mf{T}_1$ are identical to that of $\mf{T}_0$, so we may identify $\sigma'$ with a cycle of $C_k^{(0)}(\mf{T}^1)$. The homology class of this cycle is precisely $f([\sigma])$, and does not depend on the representative of $[\sigma]$ or $\tau \in C_{k+1}^{(1)}(\mf{T}^1)$ that we chose. So to define $f([\sigma])$, we may `push' $\sigma$ to the $k$-skeleton of $\mf{T}_1$ in a way which is locally determined in $\mf{T}_1$, and identify the result with the analogous homology class from $H_k^{(0)}(\mf{T}^1)$.

\begin{exmp} Recall from Example \ref{ex: 1d approximants} that in the one-dimensional case we may identify the generators of the degree zero approximant group with the admissible two-letter words of $\mc{A}^2_\omega$, and in degree one with the letters of $\mc{A}$. In degree one, as is always the case in top degree, we have that $H_1(\mf{T}^1) \cong \Z$ is generated by a fundamental class. To compute the connecting map in degree zero, let $ab$ be an admissible two-letter word, which represents an indicator chain of $C_0^{(0)}(\mf{T}^1)$ (which, abusing notation, we shall also name $ab$ here). Considered as a chain of the supertiling, $ab$ lifts to the element $\iota(ab)$ which marks each $xy$ vertex of the supertiling for which the last letter of $x$ substitutes to $a$ and the first letter of $b$ substitutes to $y$, as well as interior vertices of the supertiles, corresponding, for a supertile with label $x$, to occurrences of $ab$ in $\omega(x)$. There exists a $1$-chain $\tau$ of the original CW decomposition of the tiling which only depends on ambient supertiles and for which $ab + \partial(\tau)$ is supported on the $0$-skeleton of the supertiling. For example, we may choose $\tau$ so as to shift all $ab$ vertices of $\mf{T}_0$ contained in the interiors of supertiles to the right endpoints of these supertiles. So $f([ab])$ is represented by the chain
\[
\sigma \coloneqq \sum_{xy \in \mc{A}_\omega^2} (\omega_\text{left}^{ab}(xy)) \cdot xy
\]
where $\omega_\text{left}^{ab}(xy)$ is the number of occurrences of the word $ab$ in the substituted word $\omega(x) \cdot \omega(y)$ with first letter of that occurrence of $ab$ lying to the left of the $\cdot$ place-holder. In the notation of the definition of the connecting map $f$, we have that $q(\sigma) = \iota(ab) + \partial(\tau)$, so $f([ab]) \coloneqq q_*^{-1}(\iota_*([ab])) = [\sigma]$. Since those cycles associated to indicator cochains of admissible two-letter words generate $H_0^{(0)}(\mf{T}^1)$, the above rule determines the connecting map $f$.\end{exmp}

\begin{exmp} The approximant homology groups for the Penrose kite and dart tilings are free Abelian of ranks $8$, $5$ and $1$ in degrees $0$, $1$ and $2$, respectively. The connecting map turns out to be an isomorphism in each degree. There is a subtlety with the Penrose kite and dart substitution \cite{Pen84}, in that the substituted tiles have larger support than the inflated prototiles. In particular, the cell complex $\mc{T}_0$ of the tiling does not refine the complex $\mc{T}_1$ of the supertiling. This is only a minor inconvenience; one may work over a finer complex, corresponding to a Robinson triangle tiling, which refines both. The general procedure described above remains essentially the same, and we shall subdue this point in our discussion.

To demonstrate a typical application of the connecting map, we shall consider how it acts on the cycle $\rho' \in C_1^{(0)}(\mf{T}^1)$ of Example \ref{ex: Penrose PE}, illustrated in red in Figure \ref{fig: Penrose}, which trails the bottoms of the dart tiles. We firstly consider $\rho'$ as a $1$-cycle of the complex $\mc{T}_0$ which only depends at any given $1$-cell by those \emph{supertiles} of $\mf{T}_1$ properly containing it; formally we consider the chain $\iota(\rho') \in C_1^{(1)}(\mf{T}^1)$. Let $\tau \in C_2^{(1)}(\mf{T}^1)$ be the indicator $2$-chain of the dart tiles of $\mf{T}_0$; it is the blue chain of Figure \ref{fig: Penrose} (of course, we in fact have that $\tau$ is a member of the subcomplex $C_2^{(0)}(\mf{T}^1)$, that is, $\tau$ only depends on ambient \emph{tiles}, rather than \emph{super}tiles in this case). Then $\rho' + \partial(\tau)$ is the $1$-cycle, illustrated in green in Figure \ref{fig: Penrose}, which runs along the $1$-cells at the bottoms of the super-dart tiles, but with the opposite corresponding orientation to $\rho'$. That is to say, $\iota(\rho') + \partial(\tau) = q(-\rho')$, so $f([\rho']) \coloneqq q_*^{-1}(\iota_*([\rho'])) = -[\rho']$. More informally, we `push' the $1$-cycle $\rho'$ to the $1$-skeleton of the supertiling by adding to it the boundary of $\tau$, which is defined at any $2$-cell by only which supertiles contain it, and identify the result with the corresponding homology class of $H_1^{(0)}(\mf{T}^1)$.\end{exmp}

\subsection{Generalisations} \label{subsec: Generalisations}  There are several ways in which the method discussed above may be generalised, and conditions of the substitution rule which may easily be relaxed. For example, the primtivity condition of the substitution and the compatibility of the substitution with the cellular decomposition may be weakened. More significantly, the method may be modified to apply to mixed substitution systems, to compute the ePE homology groups and applies naturally to non-Euclidean hierarchical tilings. Instead of providing the full details of each generalisation, which is not our main focus here, we shall mostly give brief outlines of the changes that need to be made in each case; the adaptations needed to the proofs of the analogues of Theorem \ref{thm: method} are relatively straightforward in each case.

\subsubsection{Mixed Substitutions} A \emph{mixed} or \emph{multi substitution system} \cite{GahMal13} is a family of substitutions acting upon the same prototile set. Loosely, whereas the language for admissible tilings of a substitution rule $\omega$ is given by iteratively applying $\omega$ to the prototiles, in a mixed system one builds the language by applying the family of substitutions to the prototiles in some chosen sequence.

Passing to the setting of mixed substitutions adds far more generality. For example, the Sturmian words associated to some irrational number $\alpha$ may be defined using a mixed substitution system, which will be purely substitutive if and only if $\alpha$ is a quadratic irrational. In contrast to the purely substitutive case, the family of one-dimensional mixed substitution tilings exhibit an uncountable number of distinct isomorphism classes of degree one \v{C}ech cohomology groups \cite{Rus15}. In a mixed substitution tiling the tiles group together to the supertiles, those into higher order supertiles, and so on, just as in the purely substitutive case, but now the rules connecting the various levels of the hierarchy are not constant. A general framework which captures this idea is laid out in \cite{PrFSad14}.

Passing to mixed substitution systems adds some complications, since the local combinatorics of the tilings $\mf{T}_n$ and the passage between them vary in $n$. However, the method as described easily generalises to such examples. Now one needs to compute the list of stars for each $\mf{T}_n$ to find the approximant homology groups, and the connecting maps may vary at each level. 

\begin{exmp}[Arnoux-Rauzy Sequences] \label{ex: PE AR} The Arnoux--Rauzy words were originally introduced in \cite{ArnRau91} as a generalisation of Sturmian words. Let $k \in \mathbb{N}_{\geq 2}$. The Arnoux--Rauzy substitutions are defined over the alphabet $\mc{A}_k = \{1,2,\ldots,k\}$ and the $k$ different substitutions $\rho_i$ are given by $\rho_i(j) = ji$ for $i \neq j$ and $\rho_i(i)=i$. Fix an infinite sequence $(n_i)_i = (n_0,n_1,\ldots) \in \mc{A}_k^{\mathbb{N}_0}$ for which each element of $\mc{A}_k$ occurs infinitely often. Then there exist bi-infinite \emph{Arnoux--Rauzy words} for which every finite sub-word is contained in some translate of a `supertile' $\rho_{n_0} \circ \rho_{n_1} \circ \cdots \circ \rho_{n_l}(i)$. We may consider such a word as defining a tiling of labelled unit intervals of $\R^1$. The system is recognisable, so for such a tiling $\mf{T}_0$, one may uniquely group the tiles to a tiling $\mf{T}_1$ of tiles of labelled intervals for which the substitution $\rho_{n_0}$ decomposes $\mf{T}_1$ to $\mf{T}_0$. The process may be repeated, leading to an infinite hierarchy of tilings $\mf{T}_n$ for which the substitution $\rho_{n_i}$ subdivides $\mf{T}_{i+1}$ to $\mf{T}_i$; the supertiles of these tilings become arbitrarily large as one passes up the hierarchy.

The two-letter words of $\mf{T}_i$ are the elements of $\mc{A}_k^2$ with at least one occurrence of $n_i \in \mc{A}_k$. So the approximant homology at level $i$, based upon stars of supertiles of $\mf{T}_i$, is isomorphic to $\Z^k$, freely generated by the indicator $0$--chains of vertices of the form $n_i \cdot j$ where $j \in \mc{A}_k$ is arbitrary. A simple calculation shows that, with this choice of basis, the connecting map between level $i$ and level $i+1$ is the unimodular matrix $M_i$ given by the identity matrix but with a column of $1$'s down the $n_i$th column which, incidentally, is the incidence matrix of the substitution $\rho_{n_i}$.

So the degree one \v{C}ech cohomology of the tiling space of the Arnoux--Rauzy words associated to any given sequence $(n_i)_i \in \mc{A}_k^{\mathbb{N}_0}$ is
\[\check{H}^1(\Omega^1_\mf{T}) \cong H^1(\mf{T}^1) \cong H_0(\mf{T}^1) \cong \varinjlim(\Z^k \xrightarrow{M_0} \Z^k \xrightarrow{M_1} \Z^k \xrightarrow{M_2} \cdots) \cong \Z^k.\]
It is interesting to note that the matrices above are related to continued fraction algorithms. For the $k=2$ case, the Arnoux--Rauzy words are precisely the Sturmian words. To an irrational $\alpha$, the sequence $(n_i)_i$ is chosen according to the continued fraction algorithm for $\alpha$ (see \cite[\S3.2]{Dur00}) and the sequence of matrices $M_i$ of the above direct limit determine the partial quotients of $\alpha$. Whilst the isomorphism classes of the first \v{C}ech cohomology groups do not distinguish these tiling spaces, their order structure \cite{OrmSad02} is a rich invariant. Although we shall not give full details here, features such as the order structure of the cohomology groups are preserved by the method calculation described above, via Poincar\'{e} duality.\end{exmp}

\subsubsection{Euclidean Pattern-Equivariance} The method is easily adjusted to compute the ePE homology groups, based on chains which are determined cell-wise by their $i$-coronas, for sufficiently large $i$, up to \emph{rigid motion} rather than just up to translation. The method proceeds as before, but where one took translation classes of stars one simply now takes stars up to rigid equivalence. 

Non-trivial isotropy may cause issues in this setting. If any of the stars have self-symmetries which act non-trivially on the central cells, then one may only compute over suitably divisible coefficients.

Computations of this sort are of particular interest since they provide different results to the ePE cohomology calculations. However, it is possible to modify the method to compute the ePE cohomology (and thus the \v{C}ech cohomology of the space $\Omega_\mf{T}^0$) using this method for two-dimensional tilings. The ideas follow naturally from the discussions of Subsection \ref{subsec: Restoring PD}. One computes $H_\bullet^\dagger(\mf{T}^0)$ by using a similar method but restricting to indicator chains on $0$-stars which assign coefficients which are divisible by $n$ to $0$-stars with $n$-fold rotational symmetry.

\begin{exmp} \label{ex: Penrose ePE Calculation} Instead of providing more detail on the general method, we now demonstrate it on the Penrose kite and dart tilings, which should provide sufficient detail to the general method. To begin calculation, we must firstly enumerate the list of stars up to \emph{rigid motion}, as we have already done in Figure \ref{fig: VaETypes}. Since there are seven $0$-stars, seven $1$-stars and two $2$-stars (the kites and darts) up to rigid motion, the approximant chain complex for the ePE homology is given by
\[
C_\bullet^{(0)}(\mf{T}^0) = 0 \leftarrow \Z^7 \xleftarrow{\partial_1} \Z^7 \xleftarrow{\partial_2} \Z^2 \leftarrow 0.\]
Again, the boundary maps are induced from the standard cellular boundary maps after identifying the generators of these chain groups with indicator chains of the tiling. For example,
\[
\partial_1(\cf(\text{E}1)) = \cf(\text{sun}) +\cf(\text{jack}) - \cf(\text{ace}) - \cf(\text{queen}).
\]
We calculate the approximant homology groups as
\[
H_k^{(0)}(\mf{T}^0) \cong \begin{cases}
	\Z^2 \oplus (\Z/5\Z) \text{ for } k=0;\\
	\Z \text{ for } k=1;\\
	\Z \text{ for } k=2.\\
	\end{cases}
\]

The connecting map is defined essentially identically to the translational case, and we find that it is an isomorphism in each degree. So the approximant homology groups above are isomorphic to the ePE homology groups of the Penrose kite and dart tilings, and the generators of the approximant homology groups may be taken as generators of the ePE homology.

The most interesting feature here is the $5$-torsion in degree zero, which is not found in the degree two ePE cohomology, breaking Poincar\'{e} duality over integral coefficients. It is generated by the element $t= \cf(\text{sun}) + \cf(\text{star}) - \cf(\text{queen})$, illustrated in Figure \ref{fig: S0TorBdy} where one can see that $5t$ is nullhomologous via the boundary of $-\cf(\text{E}1)+\cf(\text{E}2) - \cf(\text{E}4) -2 \cdot \cf(\text{E}7)$.

To calculate the ePE cohomology, or equivalently (according to Theorem \ref{thm: Cech = ePE}) the \v{C}ech cohomology $\check{H}^\bullet(\Omega_\mf{T}^0)$, we may compute the Poincar\'{e} dual groups $H_{2-\bullet}^\dagger(\mf{T}^0)$ (see Subsection \ref{subsec: Restoring PD}). The method is similar to before, but now one replaces the approximant complex with the subcomplex
\[
0 \leftarrow 5\Z \oplus 5\Z \oplus \Z^5 \xleftarrow{\partial_1} \Z^7 \xleftarrow{\partial_2} \Z^2 \leftarrow 0,
\]
where the degree zero chain group is the subgroup of $C_0^{(0)}(\mf{T}^0)$ which restricts the coefficients on the sun and star vertices to multiples of $5$, since these vertices have $5$--fold rotational symmetry and the other vertices have trivial rotational symmetry. We calculate the modified approximant homology groups in degree zero as $\Z^2$ and the connecting maps as isomorphisms, in agreement with the following chain of isomorphisms:
\[
\check{H}^2(\Omega^0_\mf{T}) \cong H^2(\mf{T}^0) \cong H_0^\dagger(\mf{T}^0) \cong \varinjlim(\Z^2,f) \cong \Z^2.
\]
\end{exmp}

\subsubsection{Non-Euclidean Tilings} The final generalisation which we shall discuss is to non-Euclidean tilings. Again, we shall use the pentagonal tilings of Bowers and Stephenson as our running example. There is no natural action of translation for these tilings but, as discussed in Example \ref{ex: BSP}, there is a natural analogue of the ePE (co)homology groups. The method above, with essentially no modifications, may be used to compute them:

\begin{exmp} \label{ex: ePE BSP} One begins by listing the rigid equivalence classes of stars from the tiling. For a Bowers--Stephenson pentagonal tiling, there are two $0$-stars, corresponding to those vertices meeting three tiles and those meeting four, there is one $1$-star and one $2$-star. However, the $1$-star has rotational symmetry which reverses the orientation of its central $1$-cell. Since we have non-trivial cell isotropy, we may only compute over suitably divisible coefficients. Over $\Q$ coefficients, the approximant complex is given by
\[
C_\bullet^{(0)}(\mf{T}^0) = 0 \leftarrow \Q^2 \xleftarrow{\partial_1} 0 \xleftarrow{\partial_2} \Q \leftarrow 0.
\]
It follows that the approximant homology over $\Q$ coefficients is $H_k^{(0)}(\mf{T}^0;\Q) \cong \Q^2$, $0$, $\Q$ for $k=0$, $1$, $2$, respectively. The connecting map has the analogous definition to the Euclidean case, and we find it to be an isomorphism in each degree.

To compute homology over integral coefficients, we pass to the barycentric subdivision. Now we have four $0$-stars: two of them corresponding to the two $0$-stars of the original tiling, one corresponding to the barycentre of each edge and one corresponding to the barycentre of each pentagon. There are three $1$-stars and two $2$-stars. So the approximant chain complexes over $\Z$ coefficients are
\[
C_\bullet^{(0)}(\mf{T}_\Delta^0) = 0 \leftarrow \Z^4 \xleftarrow{\partial_1} \Z^3 \xleftarrow{\partial_2} \Z^2 \leftarrow 0.
\]
One computes that $H_0^{(0)}(\mf{T}^0_\Delta) \cong \Z^2$, $H_1^{(0)}(\mf{T}^0_\Delta) \cong 0$ and $H_2^{(0)}(\mf{T}_\Delta) \cong \Z$. So $H_1(\mf{T}^0_\Delta) \cong 0$, and of course $H_2(\mf{T}^0_\Delta) \cong \Z$ is generated by a fundamental class. One may calculate the connecting map in degree zero as having eigenvectors which span $\Z^2$ and have eigenvalues $1$ and $6$, so $H_0(\mf{T}^0_\Delta) \cong \Z \oplus \Z[1/6]$.

To compute the analogue of the ePE cohomology of $\mf{T}_\Delta$ (and hence the \v{C}ech cohomology of the associated tiling space $\Omega_\mf{T}^0$), one may calculate the modified groups $H_\bullet^{\dagger}(\mf{T}_\Delta^0)$ and implement Poincar\'{e} duality. At the approximant stage, this amounts to using instead the chain complexes
\[
0 \leftarrow 2\Z \oplus 3\Z \oplus 4\Z \oplus 5\Z \xleftarrow{\partial_1} \Z^3 \xleftarrow{\partial_2} \Z^2 \leftarrow 0
\]
since the $0$-stars of $\mf{T}_\Delta$ possess isotropy of orders $2$, $3$, $4$ and $5$. After computing the connecting maps and corresponding direct limits, we find that
\[
\begin{array}{l}
\check{H}^0(\Omega^0_\mf{T}) \cong H^0(\mf{T}^0_\Delta) \cong H_2^\dagger(\mf{T}^0_\Delta) \cong \Z;\\
\check{H}^1(\Omega^0_\mf{T}) \cong H^1(\mf{T}^0_\Delta) \cong H_1^\dagger(\mf{T}^0_\Delta) \cong 0;\\
\check{H}^2(\Omega^0_\mf{T}) \cong H^2(\mf{T}^0_\Delta) \cong H_0^\dagger(\mf{T}^0_\Delta) \cong \Z \oplus \Z[1/6].
\end{array}
\]
It should be remarked that the above calculation may be performed quite painlessly by hand, the combinatorial information required being surprisingly manageable, in spite of the substitution rule not forcing the border. \end{exmp}

\subsection{Proof of Theorem \ref{thm: method}.}

\begin{definition} For a cellular Borel--Moore $k$-chain $\sigma \in C_k^{\text{BM}}(\mc{T}_0)$, write $\sigma \in C_k^{(n)}(\mf{T}^1)$ to mean that $\sigma$ is determined at any $k$-cell $c$ of $\mc{T}_0$ by the immediate surroundings of $c$ in the level $n$ supertiling $\mf{T}_n$. More precisely, whenever there is a translation mapping $k$-cell $a$ to $b$ in $\mc{T}_0$, and also the patch of tiles of $\mf{T}_n$ containing $a$ to the corresponding patch at $b$, then the $k$-cells $a$ and $b$ have the same coefficient in $\sigma$. We say that $\sigma$ is \emph{hierarchical} if $\sigma \in C_k^{(n)}(\mf{T}^1)$ for some $n \in \mathbb{N}_0$ and write $C_k^{(\infty)}(\mf{T}^1)$ to denote the collection of all hierarchical $k$-chains. \end{definition}

It is not hard to see that for $\sigma \in C_k^{(n)}(\mf{T}^1)$ we also have that $\partial(\sigma) \in C_{k-1}^{(n)}(\mf{T}^1)$, so we have chain complexes $C_\bullet^{(n)}(\mf{T}^1)$ for every $n \in \mathbb{N}_0 \cup \{\infty\}$. Furthermore, since the level $n$ supertiling $\mf{T}_n$ is determined locally by the level $(n+1)$ supertiling $\mf{T}_{n+1}$ via the substitution rule, for all $m \leq n$ we have an inclusion of chain complexes
\[
\iota_{m,n} \co C_\bullet^{(m)}(\mf{T}^1) \hookrightarrow C_\bullet^{(n)}(\mf{T}^1).
\]

Importantly, every hierarchical chain is PE. Indeed, if a chain $\sigma$ is determined locally at a cell $c$ depending only on where $c$ sits in the patch of supertiles containing $c$, then $\sigma$ is also determined there by the $i$-corona of $c$ in $\mf{T}_0$, where $i$ is chosen large enough so as to deduce the level $n$ supertile decomposition at $c$ (such an $i$ exists by recognisability). We denote the inclusion of chain complexes by
\[
\iota_\infty \co C_\bullet^{(\infty)}(\mf{T}^1) \hookrightarrow C_\bullet(\mf{T}^1).
\]
Unfortunately, it is not true that a PE chain must be hierarchical. Indeed, the $i$-corona of a cell $c$ of $\mf{T}_0$ need not be determined by the supertile containing $c$ when $c$ is interior, but close to the boundary of a supertile. On the other hand, if a $k$-cycle $\sigma$ is PE, intuitively $\sigma$ only depends on very local combinatorics of the tiles of $\mf{T}_n$, relative to the sizes of the tiles of $\mf{T}_n$, for $n$ sufficiently large. One would expect that such a $k$-cycle could be perturbed, in a pattern-equivariant way, to a hierarchical $k$-cycle. More precisely, one would expect for there to exist a PE $(k+1)$-chain $\tau$ for which $\sigma + \partial(\tau)$ is supported on the $k$-skeleton and is still PE to a small radius relative to the sizes of the tiles, which would force $\sigma + \partial(\tau)$ to be hierarchical. To this end, we introduce the following technical lemma:

\begin{lemma} Let $P$ be a $d$-dimensional polytope, with polytopal decomposition $\partial \mc{P}$ of its boundary and $\delta > 0$; there exist a function $h \colon \R_{>0} \rightarrow \R_{>0}$ satisfying the following. Let $\lambda > 1$ and the inflated polytope $\lambda P$ have a cellular decomposition $\mc{P}_\lambda$ whose cells have diameter bounded by $\delta$. For any relative $k$-cycle $\sigma$ of $\mc{P}_\lambda$ modulo the boundary $\partial \mc{P}_\lambda$, with $k<d$, there exists $\tau_\sigma \in C_{k+1}(\mc{P}_\lambda)$ for which $\sigma + \partial(\tau_\sigma)$ is supported on $\partial \mc{P}_\lambda$; we may choose such chains $\tau_\sigma$ so that, whenever $\sigma_1$ and $\sigma_2$ agree at distance greater than $r$ from a subset of cells of $\lambda(\partial\mc{P})$, then $\tau_{\sigma_1}$ and $\tau_{\sigma_2}$ agree at distance greater than $h(r)$ from those cells.\end{lemma}

\begin{proof} Let $\sigma \in C_k(\mc{P}_\lambda)$ be a chain with $\partial(\sigma)$ supported on $\partial \mc{P}_\lambda$. By the homological properties of cells, there exists some $\tau_\sigma$ for which $\sigma + \partial(\tau_\sigma)$ is supported on $\partial \mc{P}_\lambda$. For any other choice of $\sigma'$ which agrees with $\sigma$ further than distance $r$ from $\partial \mc{P}_\lambda$, we have that $\sigma' + \partial(\tau_\sigma)$ is supported on an $r+\delta$ neighbourhood of $\partial \mc{P}_\lambda$. So we may restrict attention to those relative cycles $\sigma$ supported on an $r+\delta$ neighbourhood of $\partial \mc{P}_\lambda$.

So now suppose that $\sigma$ is supported on an $r+\delta$ neighbourhood of $\partial \mc{P}_\lambda$, and let $c \in \partial \mc{P}$ be a $(d-1)$-cell of the boundary of $\mc{P}$. There exists a chain $\tau$ of $\mc{P}_\lambda$, supported on a $C r + \delta$ neighbourhood of $\lambda c$, for which $\sigma + \partial(\tau)$ is supported on $\partial \mc{P}_\lambda$ union a $Cr + \delta$ neighbourhood of the remaining $(d-1)$ cells of $\lambda(\partial\mc{P})$; here, $C$ only depends on the polytope $\mc{P}$. Indeed, for sufficiently large $\lambda$, this statement would hold for $\mc{P}_\lambda$ replaced with a cellular decomposition of the disc of radius $\lambda$, the chain $\tau$ being induced by a radial deformation retraction of an $r+\delta$ neighbourhood of the boundary of the disc to its boundary. The result for the polytope $P$ may be lifted from the case of a disc by the fact that polytopes are bi-Lipschitz equivalent to the standard unit disc. We may repeat this construction for the remaining $(d-1)$-cells. As a result, we construct a chain $\tau_\sigma$ for which $\sigma + \partial(\tau_\sigma)$ is supported on $\partial \mc{P}_\lambda$. For any relative cycle $\sigma'$ which agrees with $\sigma$ distance further than $r_1$ away from the $(d-2)$-skeleton of $\lambda \partial \mc{P}$, we have that $\sigma' + \partial(\tau_\sigma)$ is supported on some $r_2$ neighbourhood of the $(d-2)$-skeleton, where $r_1$ and $r_2$ depend only on $r$ and not $\lambda$. We may now repeat this argument for those relative chains which are supported on neighbourhoods (of radius depending only on $r$) of the $k$-skeleton of $\lambda(\partial \mc{P})$ for successively smaller $k$, from which the result follows. \end{proof}

\begin{lemma} \label{lem: hier QI} The inclusion $\iota_\infty \co C^{(\infty)}_\bullet(\mf{T}^1) \hookrightarrow C_\bullet(\mf{T}^1)$ is a quasi-isomorphism.\end{lemma}

\begin{proof} For a cellular Borel--Moore $k$-chain $\sigma$ of $\mc{T}_0$, let us write that $\sigma$ is $\text{PE}_n(r)$ to mean that the values of $\sigma$ at two $k$-cells $a$ and $b$ of $\mc{T}_0$ are equal whenever there are points $x \in a$ and $y \in b$ (as open cells) for which the patches $\mf{T}_n[B_r + x]$ and $\mf{T}_n[B_r + y]$ are equal up to a translation taking $a$ to $b$. So a hierarchical chain is nothing other than a chain which is $\text{PE}_n(0)$ for some $n \in \N_0$.

Let $\sigma \in C_k(\mf{T}^1)$, so $\sigma$ is $\text{PE}_n(r)$ for all $n \in \N_0$ for some $r$. By the fact that the cells are polytopal, for sufficiently large $n$ we have that the patch of tiles of $\mf{T}_n$ within radius $r$ of any cell $c$ of $\mc{T}_0$ is determined by the star in $\mf{T}_n$ of a cell of $\mc{T}_n$ within $Cr$ of $c$. So the value of $\sigma$ at any cell $c$ of $\mc{T}_0$ contained in a supertile $t$ of $\mf{T}_n$ is determined by the star of a subcell of $t$ in $\mf{T}_n$ within radius $Cr$ of $c$.

We may now appeal to the previous lemma. By FLC, there are only a finite number of translation classes of polytopal cells in $\mf{T}_0$, and they all have diameter bounded by some $\delta > 0$. We may find a function $h \colon \R_{>0} \rightarrow \R_{>0}$ satisfying the result of the above lemma for each polytopal cell. Since $\sigma$ only depends at subcells of supertiles on stars of cells within $Cr$ of those cells, we may construct a chain $\tau$ for which $\sigma + \partial(\tau)$ is supported on the $(d-1)$-skeleton of $\mf{T}_n$ and which only depends on stars of cells within $h(r)$ in $\mf{T}_n$. So we may find a chain $\tau$ for which $\sigma + \partial(\tau)$ is supported on the $(d-1)$-skeleton of $\mf{T}_n$ and which is $\text{PE}_n(r')$, where $r'$ only depends on $r$ and not on $n$.

For sufficiently large $n$, we may repeat this process down the skeleta. As a result, we construct a chain $\tau$ for which $\sigma + \partial(\tau)$ is supported on the $k$-skeleton of $\mf{T}_n$ (recall that $k$ is the degree of $\sigma$) and is $\text{PE}_n(r')$ with $r'$ depending only on $r$. We claim that for sufficiently large $n$ such a cycle must be hierarchical. Indeed, since $\sigma$ is a cycle, it is determined across any $k$-cell of $\mf{T}_n$ by its value on any $k$-cell of $\mc{T}_0$ contained in that cell. For sufficiently large $n$, for each $k$-cell $c$ of $\mc{T}_n$ there is an interior $k$-cell of $\mc{T}_0$ for which the tiles of $\mf{T}_n$ within radius $r'$ about that cell are precisely the star of $c$ in $\mf{T}_n$. Since $\sigma + \partial(\tau)$ is $\text{PE}_n(r')$, it follows that $\sigma$ is determined at any $k$-cell of $\mf{T}_n$ by the star of that cell. Hence, $\sigma + \partial(\tau) \in C_k^{(\infty)}(\mf{T}^1)$, so we have shown that $\iota_*$ is surjective. Showing injectivity is analogous, applying the same procedure to boundaries in place of cycles. \end{proof}

The lemma above allows us to work with the chain complex $C_\bullet^{(\infty)}(\mf{T}^1)$ in computing the homology of $C_\bullet(\mf{T}^1)$. The advantage to this is that $C_\bullet^{(\infty)}(\mf{T}^1)$ possesses a natural filtration by the sub-chain complexes $C_\bullet^{(n)}(\mf{T}^1)$. The following lemma shows that the \emph{approximant homology groups} $H_\bullet^{(n)}(\mf{T}^1)$ of these sub-complexes are all isomorphic, and in a way such that the induced inclusions $(\iota_{n,n+1})_*$ between them are the same:

\begin{lemma} \label{lem: hierarchical refinement} For $n \in \mathbb{N}_0$ we have canonical quasi-isomorphisms
\[
q_n \co C_\bullet^{(0)}(\mf{T}^1) \hookrightarrow C_\bullet^{(n)}(\mf{T}^1)
\]
for which $(q_1)_*^{-1} \circ (\iota_{0,1})_* = (q_{n+1})_*^{-1} \circ (\iota_{n,n+1})_* \circ (q_n)_*$.
\end{lemma}

\begin{proof} Let $m \leq n$. Up to rescaling, the set of stars of $\mf{T}_m$ and $\mf{T}_n$ are identical. A generating element of $C_k^{(m)}(\mf{T}^1)$ is an indicator chain $\cf(c)$ of a translation class of $k$-cell $c$ of $\mc{T}_0$, labelled by where it lies in the patch of tiles of $\mf{T}_m$ containing it. We may define chain maps $q_{m,n} \colon C_\bullet^{(m)}(\mf{T}^1) \rightarrow C_\bullet^{(n)}(\mf{T}^1)$ by sending such an element to the sum of indicator chains $\cf(c')$ of translation classes of $k$-cells $c'$, suitably oriented, which are contained in the regions occupied by the analogous (inflated) locations of $c$ in the tiling $\mf{T}_n$.

Let us write $q_n$ for $q_{0,n}$. A $k$-cycle $\sigma \in C_k^{(n)}(\mf{T}^1)$ is in the image of $q_n$ if and only if it is supported on the $k$-skeleton of $\mf{T}_n$. We may now mimic the proof of Lemma \ref{lem: refinement}, replacing $\iota$ with $q_n$, to show that each $q_n$ is a quasi-isomorphism. The following identities are easily verified: $q_{j,k} \circ q_{i,j} = q_{i,k}$ and $\iota_{j,j+1} \circ q_j = q_{1,j+1} \circ \iota_{0,1}$ for $i \leq j \leq k$. It follows that $(q_{n+1})_*^{-1} \circ (\iota_{n,n+1})_* \circ (q_n)_* = (q_{n+1})^{-1}_* \circ (q_{1,n+1})_* \circ (\iota_{0,1})_* = (q_{n+1})_*^{-1} \circ ((q_{0,n+1})_* \circ (q_{0,1})_*^{-1}) \circ (\iota_{0,1})_* = (q_1)_*^{-1} \circ (\iota_{0,1})_*$.\end{proof}

The chain maps $\iota_{0,1}$ and $q_1$ are denoted $\iota$ and $q$, respectively, in the definition of the connecting map $f \coloneqq (q_*)^{-1} \circ \iota_*$. We may now prove Theorem \ref{thm: method}, that $H_\bullet(\mf{T}^1)$ is canonically isomorphic to the direct limit $\varinjlim (H_\bullet^{(0)}(\mf{T}^1), f)$. By Lemma \ref{lem: hierarchical refinement}, we have the following diagram
\begin{center}
\begin{tikzpicture}[auto]
\node (a1) {$H^{(0)}_\bullet(\mf{T}^1)$};
\node (a2) [right= of a1] {$H^{(0)}_\bullet(\mf{T}^1)$};
\node (a3) [right= of a2] {$H^{(0)}_\bullet(\mf{T}^1)$};
\node (a4) [right= of a3] {$H^{(0)}_\bullet(\mf{T}^1)$};
\node (a5) [right= of a4] {$H^{(0)}_\bullet(\mf{T}^1)$};
\node (a6) [right= of a5] {$\cdots$};

\node (b1) [below= of a1] {$H^{(0)}_\bullet(\mf{T}^1)$};
\node (b2) [right= of b1] {$H^{(1)}_\bullet(\mf{T}^1)$};
\node (b3) [right= of b2] {$H^{(2)}_\bullet(\mf{T}^1)$};
\node (b4) [right= of b3] {$H^{(3)}_\bullet(\mf{T}^1)$};
\node (b5) [right= of b4] {$H^{(4)}_\bullet(\mf{T}^1)$};
\node (b6) [right= of b5] {$\cdots$};

\begin{scope}[every node/.style={scale=.7}]
\draw [->] (a1) to node {$f$} (a2);
\draw [->] (a2) to node {$f$} (a3);
\draw [->] (a3) to node {$f$} (a4);
\draw [->] (a4) to node {$f$} (a5);
\draw [->] (a5) to node {$f$} (a6);

\draw [->] (b1) to node {$(\iota_{0,1})_*$} (b2);
\draw [->] (b2) to node {$(\iota_{1,2})_*$} (b3);
\draw [->] (b3) to node {$(\iota_{2,3})_*$} (b4);
\draw [->] (b4) to node {$(\iota_{3,4})_*$} (b5);
\draw [->] (b5) to node {$(\iota_{4,5})_*$} (b6);

\draw [->] (a1) to node {$(q_0)_* = \operatorname{id}$} (b1);
\draw      (a1) to node[swap] {$\cong$} (b1);
\draw [->] (a2) to node {$(q_1)_*$} (b2);
\draw      (a2) to node[swap] {$\cong$} (b2);
\draw [->] (a3) to node {$(q_2)_*$} (b3);
\draw      (a3) to node[swap] {$\cong$} (b3);
\draw [->] (a4) to node {$(q_3)_*$} (b4);
\draw      (a4) to node[swap] {$\cong$} (b4);
\draw [->] (a5) to node {$(q_4)_*$} (b5);
\draw      (a5) to node[swap] {$\cong$} (b5);
\end{scope}
\end{tikzpicture}
\end{center}
This isomorphism of directed systems induces an isomorphism
\[\varinjlim (H_\bullet^{(0)}(\mf{T}^1), f) \cong \varinjlim (H_\bullet^{(n)}(\mf{T}^1), (\iota_{m,n})_*).\]
Since, by definition, $C_\bullet^{(\infty)}(\mf{T}^1) = \bigcup_{n=0}^\infty C_\bullet^{(n)}(\mf{T}^1)$, we may identify $C_\bullet^{(\infty)}(\mf{T}^1)$ with the direct limit $\varinjlim (C_\bullet^{(n)}(\mf{T}^1), \iota_{m,n})$. So by Lemma \ref{lem: hier QI} we have the following string of quasi-isomorphisms
\[
\varinjlim (C_\bullet^{(n)}(\mf{T}^1), \iota_{m,n}) \xrightarrow{\cong} C_\bullet^{(\infty)}(\mf{T}^1) \xhookrightarrow{\iota_\infty} C_\bullet(\mf{T}^1).
\]
Applying homology and combining with the isomorphism of direct limits established above, we have that $\varinjlim (H_\bullet^{(0)}(\mf{T}^1),f) \cong H_\bullet(\mf{T}^1)$.

\end{document}